\documentclass{amsart}
\usepackage[OT2,T1]{fontenc}
  \DeclareSymbolFont{cyrletters}{OT2}{wncyr}{m}{n}
  \DeclareMathSymbol{\Sha}{\mathalpha}{cyrletters}{"58}
\usepackage{fancyhdr}
\usepackage{epsf}
\usepackage{times}
\usepackage{amsfonts}
\usepackage{graphicx}
\usepackage{amsfonts}
\usepackage{amsmath}
\usepackage{amscd}
\usepackage{amssymb}
\usepackage{tikz}
\usepackage{tikz-cd}
\usepackage{rotating}
\newcommand*{\isoarrow}[1]{\arrow[#1,"\rotatebox{90}{\(\sim\)}"
]}
\usepackage[utf8]{inputenc}
\usepackage[all]{xy}

\usepackage{color}

\usepackage{stmaryrd}
\usepackage{amsthm}
\usepackage{amssymb}
\usepackage{amsxtra}
\usepackage{fullpage}
\usepackage[colorlinks=true]{hyperref}
\usepackage{tikz} 
\usetikzlibrary{cd}
\tikzcdset{arrow style=math font}
\usepackage{microtype}

\newtheorem{theorem}{Theorem}
\newtheorem*{theorem*}{Theorem}
\newtheorem*{definition*}{Definition}
\newtheorem{proposition}[theorem]{Proposition}
\newtheorem{lemma}[theorem]{Lemma}
\newtheorem{corollary}[theorem]{Corollary}

\theoremstyle{definition}
\newtheorem{definition}[theorem]{Definition}

\newtheorem{remark}[theorem]{Remark}
\newtheorem*{remark*}{Remark}

\newtheorem{example}[theorem]{Example}

\newtheorem{conjecture}[theorem]{Conjecture}
\newtheorem*{conjecture*}{Conjecture}

\title{Values of zeta functions of arithmetic surfaces at $s=1$}
\author{Stephen Lichtenbaum and Niranjan Ramachandran}
\address{Stephen Lichtenbaum, Department of Mathematics, Brown University, Providence, RI 02912}
\email{stephen\_lichtenbaum@brown.edu}
\urladdr{\href{https://www.math.brown.edu/faculty/lichtenbaum.html}{Webpage}} 

\address{Niranjan Ramachandran, Department of Mathematics, University of Maryland, College Park, MD 20742 USA.}
\email{atma@umd.edu}
\urladdr{\href{http://www2.math.umd.edu/~atma/}{Webpage}}
\begin{document}
\begin{abstract}  We show that the conjecture of \cite{SL2020}  for the special value at $s=1$ of the zeta function of an arithmetic surface is equivalent to the Birch-Swinnerton-Dyer  conjecture for the Jacobian of the generic fibre. 
\end{abstract}
\maketitle

\section{Introduction and statement of results}

Let $V$ be a connected regular scheme with $V \to \textrm{Spec}~\mathbb Z$ flat projective and consider its scheme zeta function
\[\zeta(V,s) = \prod_{x \in V} \frac{1}{(1-N(x)^{-s})};\]the product is over all closed points $x$ of $V$ and $N(x)$ is the size of the finite residue field of $x$.  
 The first author \cite[Conjecture 3.1]{SL2020} has recently conjectured a formula for the special value $\zeta^*(V,r)$ as a generalized Euler characteristic (up to powers of two):
\begin{equation}\label{conj} \zeta^*(V,r) = \pm \chi(V,r).\end{equation}

The basic example is the case when $V = \textrm{Spec}~\mathcal O_F$ where $\mathcal O_F$ is the ring of integers of a number field $F$, then (\ref{conj})  recovers \cite[\S 7]{SL2020} the well known class number formula (\ref{cnf1}) for $r=0$ and $r=1$ (see \S \ref{d=1}) and, for $r\neq 0,1$,  the formula (\ref{conj}) recovers the conjecture \cite{MR0406981} of the first-named author relating $\zeta^*(V,r)$ and algebraic K-theory. The next case of interest would be that of arithmetic surfaces and  $r=1$ (considered in this paper). The compatibility of  (\ref{conj})  with the functional equation,  proved in \cite[\S 6]{SL2020} using the deep results of Bloch-Kato-Saito \cite{MR2102698}, is a strong point in support of the conjecture. Our results for arithmetic surfaces provide additional evidence for (\ref{conj}).

\subsection{History} It is a fascinating and important problem to give a description of these special values in terms of 
cohomological invariants of $V$.  The first conjectural attempt at such a description was given in 1991 by Fontaine and Perrin-Riou
 \cite{MR1265546, MR1206069}, building on previous work by Deligne, Beilinson, Bloch and Kato.  Roughly speaking, Bloch and Kato \cite{MR1086888} considered values at positive integers
while Deligne and Beilinson considered values at negative integers.  Fontaine and Perrin-Riou were the first to give a uniform approach valid for all integers.

The reader should be warned that the conjectures of Fontaine and Perrin-Riou \cite{MR1265546, MR1206069} do not in fact concern the scheme zeta-function $\zeta(V,s)$ above 
but rather the Hasse-Weil zeta-function $\zeta_{HW}(V_0,s)$, which is the alternating product of Hasse-Weil L-functions $L_{HW}(h^j(V_0), s)$ \cite{SDPP_1969-1970__11_2_A4_0}.  These L-functions $L_{HW}(h^j(V_0), s)$ only depend on the generic fiber $V_0$ of $V \to \mathrm{Spec}~ \mathbb Z$, while $\zeta(V,s)$ depends on the special fibers as well.  In fact Fontaine and Perrin-Riou give conjectures for the special values of the L-functions $L_{HW}(h^j(V_0), s)$, whose alternating product  yields conjectures for the Hasse-Weil zeta-function $\zeta_{HW}(V_0, s)= \prod_j (L_{HW}(h^j(V_0), s))^{(-1)^{j+1}}$.  These two zeta-functions agree if $V$ is smooth over Spec $\Bbb Z$ but not in general.

Fontaine and Perrin-Riou first introduce various (conjecturally) finite-dimensional vector spaces and make a conjecture (which we will call Conjecture FPR) giving the special-value of the (Hasse-Weil) zeta-function $\zeta_{HW}(V_0, s)$ up to a rational number in terms of determinants of maps between these vector spaces tensored with $\Bbb R$.   They then use spaces taken from $p$-adic Hodge theory to refine these conjectures to obtain a conjecture for the special value of $\zeta_{HW}(V_0, s)$ up to sign.

A crude description of the main conjecture made in 2017 by the first author \cite{SL2020} (which we will call Conjecture L)  would be to say that it refines Conjecture FPR by using canonical integral models
for the vector spaces used in that conjecture, avoiding the necessity for p-adic Hodge theory.  In fact the difference between the scheme zeta-function and the Hasse-Weil zeta-function forces minor changes in the vector spaces to  be considered.  These changes are similar to those made in  \cite{MR3874942}  for the same reason.  In 2016 Flach and Morin  \cite{MR3874942} also made a general special-values conjecture, which we will refer to as Conjecture FM1.  Conjecture FM1 does make use of p-adic Hodge theory, and is more closely related to Conjecture FPR than is Conjecture L.  In fact, Conjecture FM1 can be shown to be equivalent to Conjecture FPR if $V$ is smooth over Spec $\Bbb Z$ - see  \cite{MR3874942}.

In 2019, Flach and Morin \cite{FM2} made a new conjecture (referred to here as Conjecture FM2) which avoids $p$-adic Hodge theory and so is more closely modeled on  Conjecture L and less related to Conjecture FPR.  In fact it seems plausible that it is not too difficult to show that Conjecture L is equivalent to Conjecture FM2, except that Conjecture FM2 eliminates the 2-power indeterminacy.  Both Conjecture L and Conjecture FM2 have been shown to be compatible \cite{SL2020, FM2} with a form of the functional equation for the scheme zeta-function.  This has not been shown for either Conjecture FM1 or Conjecture FPR. 

Conjecture L is a bit more ad hoc than Conjecture FM2, but has the advantage that it involves much less elaborate machinery.  Of course, as previously mentioned, Conjecture FM2 is more precise,  but we hope that is  possible to remedy this by modifying Conjecture L using the Artin-Verdier topology.


 Conjecture FPR is {\it local}: for each prime $p$, a conjecture is formulated involving sophisticated $p$-adic Hodge theory and the corresponding local variety $V \times \textrm{Spec}~\mathbb Q_p$; and it is their combination (for all primes $p$) that determines the special value up to sign.  In contrast, Conjecture L (and FM2) is truly {\it global} as it produces a description of $\zeta^*(V,r)$ intrinsically in terms of invariants of $V$. For instance, in the case of an arithmetic surface $X$, where Conjecture  FPR calculates the order of the Tate-Shafarevich group $\Sha(J)$ of the Jacobian $J$ of the generic fibre $X_0$ of $X$, Conjecture L  for $\zeta^*(X,1)$ concerns its intrinsic counterpart: the Brauer group $\textrm{Br}(X)$. 
 
\subsection{The conjecture for arithmetic surfaces}
Let $S ={\rm Spec}~\mathcal O_F$ be the spectrum of the ring of integers $\mathcal O_F$ in a number field $F$.  In this paper, {\it an arithmetic surface} $X$ over $S$ will mean

\begin{itemize} 
\item a regular scheme $X$ (of dimension two) together with a projective flat morphism $\pi: X \to S$;
\item the generic fiber $X_0 \to {\rm Spec}~F$ is a geometrically connected smooth projective curve of genus $g$.
\end{itemize} 
We recall the conjecture of \cite[Conjecture 3.1]{SL2020} for the case at hand: $V=X$ and $r=1$. 
\begin{conjecture}\label{main conjecture}  For an arithmetic surface $X$ as above, the special value of $\zeta^*(X,1)$ at $s=1$ is a  generalized Euler characteristic $\chi(X,1)$ (up to  powers of two)
\[\zeta^*(X,1) =  \pm \chi(X,1) = \pm \frac{\chi_{A,C}(X,1)}{\chi_B(X,1)} .\]
\end{conjecture} 
As we shall see in \S \ref{ATooo}, this provides a description of $\zeta^*(X,1)$ in terms of (periods computed using) the finitely generated $\mathcal O_F$-modules $H^*(X, \mathcal O)$, the Betti cohomology $H^*_B(X_{\mathbb C}, \mathbb Z(1))$, the conjecturally finite group $\textrm{Br}(X)$, the finitely generated abelian group $\textrm{Pic}(X)$, and the Arakelov intersection pairing on $\textrm{Pic}^0(X)$. 

Our main result is the following theorem:
 \begin{theorem}\label{Main theorem}  Conjecture \ref{main conjecture} is equivalent to the Birch-Swinnerton-Dyer conjecture (Conjecture \ref{bsd}) for the Jacobian $J$ of $X_0$ over $F$.   
 \end{theorem} 

\begin{remark} (i)  This theorem is the analogue (\S \ref{ATooo}) for arithmetic surfaces of Conjecture (d) \cite[p.~427]{Tate1966}; see  \cite{SL1983, MR528839, MR2092767, MR3858404}. Such an arithmetic analogue had long been expected.  Flach-Siebel \cite{flachsiebel} prove a similar relation between the conjectures of Flach-Morin \cite{MR3874942} for $X$ and Conjecture \ref{bsd}.  

(ii) In a certain sense, Conjecture \ref{main conjecture}  is dual to Conjecture \ref{bsd}. When $X_0$ is an elliptic curve, then the former uses the proper map $\pi:X\to S$ that may not be smooth whereas the latter uses the N\'eron model $\mathcal N$ of $X_0$ (and the smooth map $\mathcal N \to S$ may not be proper). This suggests that one should start with the case when $\pi$ is smooth and admits a section; this case is considered in \S \ref{Proof}. 

(iii) If $g=0$, then $J$ is trivial. Since Conjecture \ref{bsd} is obviously true in this case, Theorem \ref{Main theorem} gives a proof of Conjecture \ref{main conjecture} for $g=0$. See Remark \ref{giszero} for details. \qed \end{remark} 

Since Tate \cite{Tate1966} has proved that Conjecture \ref{bsd} is invariant under isogeny, one obtains the following

\begin{corollary} Suppose $X$ and $X'$ are arithmetic surfaces over $S$ such that the Jacobians of $X'_0$ and $X_0$ are isogenous over $F$. Then Conjecture \ref{main conjecture} for $X$ is equivalent to Conjecture \ref{main conjecture} for $X'$.
\end{corollary} 
 In particular, Conjecture \ref{main conjecture} for (the regular integral models of) an elliptic curve $E$ and a torsor $T$ over $E$ are equivalent.

\subsection{Plan of the paper} 

We start in  \S \ref{conjecture-section} with a presentation of  the main conjecture from \cite{SL2020} and provide a detailed description of the basic examples ($V=S$); in particular, we show that the conjecture recovers the analytic class number formula (\ref{cnf1}).  This is followed in \S \ref{Background} by the conjecture (Conjecture \ref{bsd}) of Birch-Swinnerton-Dyer (and Tate) for the global $L$-function of an abelian variety over a number field. Here we rely on \cite{Gross82} because it contains the formulation of Conjecture \ref{bsd} in terms of N\'eron models and provides the definition of a period even when the (projective) module of invariant differentials on the N\'eron model is not free as a $\mathcal O_F$-module.

 In \S \ref{2periods}, we prove the key relation (Theorem \ref{allperiods}) between the period in Conjecture \ref{bsd} and the period  in Conjecture \ref{main conjecture}. 

\S \ref{Proof} contains the proof of Theorem \ref{Main theorem} in the special but illustrative and important case where $\pi:X \to S$ is smooth and $X_0(F)$ is non-empty. 
The rich intersection theory on a regular surface figures prominently in \S \ref{tools+}. The beautiful results of Raynaud-Liu-Lorenzini \cite{MR2092767, MR1717533, BLR} on the N\'eron model, the connected components of the special fiber and the relation to intersection theory play a key role in our paper. Though the analysis of the role of torsion in derived de Rham cohomology could complicate the study of $\zeta^*(V,r)$ in general, these results help us bypass the torsion in $H^1(X, \mathcal O_X)$. An equally important role is played by the result of Geisser \cite{Geisser1}. Finally, the proof of Theorem \ref{Main theorem} is completed in \S \ref{Proof2}.   

The concept of integral structures and derived differentials is resonant with the methods used in the study of $\epsilon$-constants, Arakelov theory \cite{MR1477762, MR1914000} and Bloch's conjectural formula for the conductor \cite{MR2102698}.

\subsection{Acknowledgements} This work was initiated at the Hausdorff Research Institute during the 2018 Trimester on Periods; we thank the organizers of the Trimester and Hausdorff Research Institute for mathematics for their support. We also thank the organizers of the Motives in Tokyo conference in 2019. The referee deserves special thanks for a detailed report containing many suggestions for improvement. We are grateful to Takashi Suzuki for his comments and suggestions. In particular, the proof in Remark \ref{suzuki} is due to him.

\subsection*{Notations} 

As usual, we write $d_F$ for the discriminant of $F$ over $\mathbb Q$, $h$ for the class number of  $F$ and $w$ the number of roots of unity in $F$. Using the notation 
$A_{tor}$ for the torsion subgroup of an abelian group $A$ and $[B]$ for the cardinality of a finite group $B$, one has $w= [F^*_{tor}]$. 

We shall use $S$ to indicate both Spec~$\mathcal O_F$ and the set of non-zero prime ideals of $\mathcal O_F$. Thus, $v\in S$ means a closed point of $S$ or  a non-zero prime ideal in $\mathcal O_F$. We write $q_v$ for the size of the finite residue field $k(v)$. 

Let $\tilde{S}$ denote the set of imbeddings of $F$ into $\mathbb C$; complex conjugation $c$ on $\mathbb C$ induces an involution $\sigma \mapsto c\sigma = \bar{\sigma}$ on $\tilde{S}$ via 
\[\sigma: F\to \mathbb C\qquad \mapsto \qquad {\bar{\sigma}}: F \xrightarrow{\sigma} \mathbb C \xrightarrow{c} \mathbb C.\]
The set $S_{\infty}$  of infinite places of $S$ is the set of orbits of the involution; the set of fixed points of the involution is exactly the subset $S_{\mathbb R}$ of imbeddings of $F$ into $\mathbb R$. As usual, $r_1$ is the size of $S_{\mathbb R}$ and $r
_1+ r_2$ is the size of $S_{\infty}$. 

For any scheme $Y$ over $S$, we put $Y_{\mathbb C} = Y \times_{\textrm{Spec}~\mathbb Z} \textrm{Spec}~\mathbb C$ and $Y_{\sigma}$ for the base change along  $\sigma:F \to \mathbb C$: 
\begin{equation}
\begin{tikzcd} 
Y_{\mathbb C} \arrow[r,] \arrow[ d] & Y \arrow[d]& &&Y_{\sigma} \arrow[d] \arrow{r} & Y  \arrow[ d] \\
\textrm{Spec}~\mathbb C \arrow[r] & \textrm{Spec}~\mathbb Z &&&\textrm{Spec}~\mathbb C \arrow[r] & S
\end{tikzcd} 
\end{equation}

We write $H^*_{et}(Y, \mathbb Z(r))$ for the \'etale motivic cohomology groups of $Y$; these are the hypercohomology groups of Bloch's higher Chow group complex in the \'etale topology of $Y$.

\section{Statement of the main conjecture from \cite{SL2020}}\label{conjecture-section}

Fix  a connected regular scheme $V$ of dimension $d$ together with a projective flat morphism $V \to S$; we assume that the generic fiber $V_0 \to \textrm{Spec}~F$ is a geometrically connected smooth variety of dimension $d-1$ over Spec~$F$.  The aim of this section is to present the conjectural formula from \cite{SL2020} for the special value $\zeta^*(V,r)$. 

\subsection{The Bloch-Kato-Fontaine-Perrin-Riou conjecture}\label{BKFP} \cite{MR1206069, MR1265546} 

 The motive $M= h^j(V_0)(r)$ has a Betti realization $H^j_B(V
 _{\mathbb C}, \mathbb Z(r))$,  a de Rham realization $H^j_{dR} (V_{\mathbb C}, \mathbb C(r))$, an action of complex conjugation $c$ on $H^j_B(V
 _{\mathbb C}, \mathbb Z(0)) = H^j_B(V_{\mathbb C}, \mathbb Z)$,  a Hodge filtration $F$ on  $H^j_{dR} (V_{\mathbb C}, \mathbb C(r))$, and a period map
\[ \theta_{j,r}:  H^j(V_{\mathbb C}, \mathbb C(r)) \xrightarrow{\sim} H^j_{dR}(V_{\mathbb C}, \mathbb C(r)).\]
As abelian groups, $H^j_B(V_{\mathbb C}, \mathbb Z(r)) = H^j_B(V_{\mathbb C}, \mathbb Z)$ and $H^j_B(V _{\mathbb C}, \mathbb C(r)) = H^j_B(V_{\mathbb C}, \mathbb C)$; also  $H^j_B(V_{\mathbb C}, \mathbb C(r)) =H^j_B(V_{\mathbb C}, \mathbb Z(r))\otimes_{\mathbb Z}\mathbb C$. The map $\theta_{j,0}$ is the classical period map from  
$H^j_B(V_{\mathbb C}, \mathbb C)$ to $H^j_{dR} (V_{\mathbb C}, \mathbb C)$. 
We define 
\[H^j(V_{\mathbb C}, \mathbb Z(r))^{+}  =  \{x \in H^j(V_{\mathbb C}, \mathbb Z)~|~ c(x) =~(-1)^r\cdot x\},\qquad H^j(V_{\mathbb C}, \mathbb C(r))^{+}  =  \{x \in H^j(V_{\mathbb C}, \mathbb C)~|~ c(x) =~(-1)^r\cdot x\}.\]

Finally, the map $\theta_{j,r}$ defined as $(2\pi i)^r\cdot\theta_{j,0}$  induces a map
\begin{equation}\label{alphaM}  \alpha_M: H^j(V_{\mathbb C}, \mathbb C(r))^{+} \to M_B \xrightarrow{\sim} M_{dR} \to t_M= \frac{M_{dR}}{F^0}.\end{equation}
Fontaine-Perrin-Riou  \cite{MR1265546, MR1206069} introduce certain vector spaces $H^0_f(M)$ and $H^1_f(M)$ that  are conjecturally  finite-dimensional $\mathbb Q$-vector spaces. 
Further, they  conjecture the existence of {\it the fundamental exact sequence} 
\begin{equation}\label{FPles}  0 \to H^0_f(M)\otimes{\mathbb R} \xrightarrow{c} \textrm{Ker}(\alpha_M) \to H^1_f(M^*(1))^*\otimes{\mathbb R}  \xrightarrow{h} H^1_f(M)\otimes{\mathbb R}  \xrightarrow{r} \textrm{Coker}(\alpha_M) \to H^0_f(M^*(1))^*\otimes{\mathbb R}  \to 0,
\end{equation}where $c$ is a cycle class map, $h$ is a (Arakelov) height pairing and $r$ is the Beilinson regulator. They used it to reformulate the Deligne-Beilinson conjecture $C_{DB}(M)$ \cite{MR1206069} which predicts the special value $L^*(M,0) \in {\mathbb R^{\times}}/{\mathbb Q^{\times}}$ at $s=0$ of the $L$-function $L(M,s)$. 

To predict $L^*(M,0) \in \mathbb R^{\times}$, they introduce similar exact sequences for each prime $p$; this is the conjecture $C_{BK}(M)$ \cite{MR1206069}  for $p$. One obtains a prediction for the special value of $\zeta_{HW}(V_0,s)$ at $s=r$ by combining $C_{BK}(M)$ for the motives $h^j(V_0)(r)$ for all $j$.  For more details, see  \cite{MR1086888, MR1265546, MR2019010, MR2882695, MR1206069, MR2088713, MR2520461, MR3874942}.

But our approach, following \cite{SL2020}, to $\zeta^*(V,r)$ is different: we endow the vector spaces in (\ref{FPles}) with integral structures using Weil-\'etale motivic cohomology groups and the de Rham cohomology groups of $V$. The determinants of the maps in (\ref{FPles}) with respect to these integral structures will be used for a description of $\zeta^*(V,r)$.

\subsection{Weil-\'etale motivic cohomology groups}These groups  $H^i_W(V, \mathbb Z(r))$ \cite[\S 2.1]{SL2020} are defined as \'etale motivic cohomology $H^i_{et}(V, \mathbb Z(r))$ for $i \le 2r$ and $r \ge 0$ and then, for $i >2r$, as the dual of \'etale motivic cohomology  
\[H^i(\textrm{RHom}(R\Gamma_{et}(V, \mathbb Z(d-r)), \mathbb Z [-2d-1])),\] thereby, for $i >2r$, they sit in an exact sequence 
\[0 \to \textrm{Ext}^1(H^{2d+2-i}_{et}(V, \mathbb Z(d-r)), \mathbb Z) \to H^i_W(V, \mathbb Z(r)) \to \textrm{Hom}(H^{2d+1-i}_{et}(V, \mathbb Z(d-r)), \mathbb Z) \to 0.\]
We define $H^i_W(X, \mathbb Z(r))$ to be zero when $r<0$ and $i\le 2r$.

The groups  $H^i_W(V, \mathbb Z(r))$ are conjectured to be finitely generated.

\subsection{Integral structures and complexes}\label{integral1} We recall the well known theory of determinants \cite[Lecture 1, \S 5]{MR2882695}. 

 A lattice $L$ in a finite dimensional vector space $C$ is the $\mathbb Z$-submodule generated by a basis of $C$.
  Write $\Lambda L \cong \mathbb Z$ (non-canonically) for the highest exterior power of $L$ and $\Lambda C \cong \mathbb C$ (non-canonically) for the highest exterior power of $C$, and $(\Lambda C)^{-1}$ for its dual.  Given a finite complex of finite dimensional complex vector spaces $C_i$ and lattices $L_i \subset C_i$
\[C_{\bullet} :\qquad 0 \to \cdots C_i \xrightarrow{f_i} C_{i+1} \xrightarrow{f_{i+1}} C_{i+2} \to \cdots, \]
one has a canonical map $g: \Lambda (L_{\bullet}) \to \Lambda (C_{\bullet})$ where
\[ \Lambda (L_{\bullet}) = \otimes_{j}(\Lambda L_j)^{(-1)^j} \cong \mathbb Z, \qquad  \Lambda(C_{\bullet}) = \otimes_{i} (\Lambda C_i)^{(-1)^i}.\]

 If $C_{\bullet}$ is exact, then one has a canonical isomorphism 
\[h: \Lambda(C_{\bullet}) \xrightarrow{\sim} \mathbb C.\]We define 
\[ det(C_{\bullet},L_{\bullet}) = h(g(e))\] as the image of a generator $e$ of $\Lambda (L_{\bullet}) \cong \mathbb Z$ under the composite map 
\[\Lambda (L_{\bullet}) \xrightarrow{g}   \Lambda (C_{\bullet}) \xrightarrow{h}  \mathbb C;\]
 $det(C_{\bullet},L_{\bullet})$ is well-defined as an element of \[ {\mathbb C^{\times}}/{\{\pm 1\}}.\] 

\subsubsection{Integral structures}  
\begin{definition} (i) An integral structure on a complex vector space $V$ is a pair $(M,m)$ consisting of a lattice $M\subset V$  and a positive rational number $m\in \mathbb Q_{>0}$. 

(ii) The Euler characteristic of a map from $(M,m)$ on $V$ to $(N,n)$ on $W$ is the determinant of the map $f: V \to W$ (with respect to the lattices $M$ and $N$) multiplied by $n/m$.

(iii) Two integral structures $(M,m)$ and $(M',m')$ on $V$ are Euler-equivalent  if the Euler characteristic (relative to $M$ and $M'$) of the identity map on $V$ is one. 

(iv) The $\mathbb Z$ - dual of  $(M,m)$  is the integral structure $(\textrm{Hom}_{\mathbb Z}(M,\mathbb Z), \frac{1}{m})$ on $\textrm{Hom}_{\mathbb C}(V, \mathbb C)$.\qed
\end{definition} 
 The Euler-characteristic of $f$ in (ii)  does not change if $(M,m)$ is replaced with 
a Euler-equivalent integral structure $(M',m')$ on $V$. If $M'$ is a sublattice of $M$ with index $e$, then $(M', em)$ and $(M,m)$ are Euler-equivalent as are $(\textrm{Hom}_{\mathbb Z}(M,\mathbb Z), 1)$  and $(\textrm{Hom}_{\mathbb Z}(M',\mathbb Z), 1/e)$: Euler-equivalence is compatible with duality.

If $\phi:A \to V$ is a homomorphism, $A$ is a finitely generated abelian group, $\textrm{Ker}(\phi)$ is finite and $\phi(A)$ is a lattice in $V$, then $(\phi(A), [\textrm{Ker}(\phi)])$ is an integral structure. If $A$ is torsion-free, $(\phi(A), 1)$ essentially determines $\phi$.     

\subsubsection{Integral structures on complexes} Given a finite complex $C_{\bullet}$ of finite dimensional complex vector spaces, an {\it integral structure} $L$ on $C_{\bullet}$ is 
a collection $L_{\bullet} = (L_i, t_i)$ where $L_i$ is a lattice in $C_i$ and $t_i$ is a positive rational number. If $C_{\bullet}$ is exact, we define 
\[\chi (C_{\bullet}, L_{\bullet}) = \frac{det (C_{\bullet}, L_{\bullet} )}{T} \in {\mathbb C^{\times}}/{\{\pm 1\}}, \qquad T = \frac{t_0\cdot ~\cdot~t_2~.\cdots}{t_1~\cdot~t_3~\cdot~\cdots} \in \mathbb Q_{>0}.\]
Here $T$ is a multiplicative Euler characteristic.

\begin{example} A collection $\phi$ of homomorphisms $\phi_i: A_i \to C_i$ where, for all $i$, $A_i$ is a finitely generated abelian group, $\textrm{Ker}(\phi_i)$  is finite, and $\textrm{Im}(\phi_i)$ is a lattice in $C_i$ gives an integral structure $L_{\bullet}=  (L_i, t_i) $ where $L_i= \phi_i(A_i)$ and $t_i =[\textrm{Ker}(\phi_i)]$; in this case, we write $T= \chi (A_{tor})$.  If $C_{\bullet}$ is exact, then we have 
\[\chi (C_{\bullet}, A_{\bullet}) = \frac{det (C_{\bullet}, L_{\bullet} )}{ \chi (A_{tor})} \in {\mathbb C^{\times}}/{\{\pm 1\}}.\]
We often refer to $A_{\bullet}$ as an integral structure on $C_{\bullet}$. 
\end{example} 

If there is a canonical integral structure $L_{\bullet}$ on $C_{\bullet}$, then we write $det(C_{\bullet})$ and $\chi(C_{\bullet})$ for $det(C_{\bullet},L_{\bullet})$ and $\chi(C_{\bullet},L_{\bullet})$.

\subsubsection{Integral structures on $\mathcal O_F$-modules} 
Let $M$ be a finitely generated  $\mathcal O_F$-module of rank $d$; there are two natural  integral structures $i_M$ and $j_M$ on $V= M\otimes_{\mathbb Z}\mathbb C $. Let $\phi: M \to V$ be the canonical map. The integral structure $i_M$ is obtained by entirely forgetting the $\mathcal O_F$-structure:
\begin{equation}\label{im-integral} 
i_M = (\phi(M), [M_{tor}]).
\end{equation} 

While $i_M$ appears in the conjectural formula (\ref{conj}), another\footnote{Recall that two lattices $L$ and $L'$ in $V$ are commensurable if $L\cap L'$ is finite index in $L$ and $L'$. The integral structures $i_M$ and $j_M$ for $M =\mathcal O_F$ are not commensurable in general.} integral structure $j_M$  appears in Conjecture \ref{bsd}; $j_M$ uses the $\mathcal O_F$-structure on $M$. Recall  the standard isomorphism
 \begin{equation}\label{standard-psi} 
\psi: \mathcal O_F \otimes_{\mathbb Z}\mathbb C \xrightarrow{\sim}  \mathbb C^{r_1+2r_2} =  \mathbb Z^{r_1+2r_2}\otimes\mathbb C
 \end{equation}  that sends $x\otimes 1\in \mathcal O_F\otimes\mathbb C$ to the element of $\mathbb C^{r_1+2r_2}$ given by the collection of $\sigma(x)$ as $\sigma\in \tilde{S}$ runs through the embeddings of $F$ in $\mathbb C$. Define $M_{\sigma} = M\otimes_{\mathcal O_F} \mathbb C$ using $\sigma: F\to \mathbb C$. 
   
The isomorphism (\ref{standard-psi}) shows
\[V= M\otimes_{\mathbb Z}\mathbb C \cong M \otimes_{\mathcal O_F} (\mathcal O_F\otimes_{\mathbb Z} \mathbb C) \cong M\otimes_{\mathcal O_F}\mathbb C^{r_1+2r_2} \cong \prod_{\sigma \in \tilde{S}} M_{\sigma}.\]

Suppose $B=\{x_1, \cdots, x_d\} \subset M$ is linearly independent over $\mathcal O_F$.  Let $\{\epsilon_j\}$ be the standard basis of $ \mathbb C^{r_1+2r_2}$. 
Define $\tilde{B} = <x_i\otimes\epsilon_j>$ to be the $\mathbb Z$-span of $x_i\otimes \epsilon_j$ inside $V$. Consider the integral structure $i(B)= (\tilde{B}, b)$ where $b$ is the order of the quotient of $M$ by the $\mathcal O_F$-submodule $<x_1, \cdots, x_d>$ generated by $B$.

\begin{proposition}\label{jm-prop}  For any two such subsets $B$ and $B'$ of $M$,  $i(B)$ and $i(B')$ are Euler-equivalent.  
\end{proposition} 
\begin{proof} Let  $B' =\{c_1x_1, c_2x_2, \cdots, c_dx_d\}$ and $C = c_1.c_2\cdots c_d$ and consider $i(B')$. As $b$ and $b'$ are related by $b' = C.b$ and the determinant of the identity map of $V$ with respect to $i(B)$ and $i(B')$ is $C$, so $i(B)$ and $i(B')$ are Euler-equivalent. Given $B$ and $B'$, one can find $B''$ such that the previous argument applies to the pairs $B, B''$ and $B', B''$.   
\end{proof} 
\begin{definition}\label{jm-integral} For any finitely generated $\mathcal O_F$-module $M$,  the  $\mathcal O_F$-integral structure $j_M$ on $V =M\otimes_{\mathbb Z}\mathbb C$
is the Euler equivalence class in Proposition \ref{jm-prop}. \qed
\end{definition}

\subsubsection{Results on integral structures} Let $M$ be a finitely generated $\mathcal O_F$-module of rank $d$. 

\begin{proposition}\label{euler-disc} (i) The Euler characteristic of the identity map on $V$ from  $i_M$ to $j_M$ is $(\sqrt{d_F})^d$. 

(ii) Let  $M$ be torsion-free and $M^{\vee}= \textrm{Hom}_{\mathcal O_F}(M, \mathcal O_F)$. The $\mathbb Z$-dual of $j_M$ is the $\mathcal O_F$-integral structure $j_{M^{\vee}}$.
\end{proposition} 

\begin{proof} (i) Write $\kappa_M$ for the Euler characteristic. For $M =\mathcal O_F$, this says $\sqrt{d_F}$ is the determinant of the standard map $\psi: \mathcal O_F\otimes_{\mathbb Z} \mathbb C \to \mathbb Z^{r_1+2r_2}\otimes\mathbb C$  in (\ref{standard-psi}). So $\kappa_M= \sqrt{d_F}^d$ if $M$ is free of finite rank $d>0$. Even if $M$ is not free, it contains a free submodule $M'$ of finite-index.  As $\kappa_M = \kappa_{M'}$, the claim for $M'$ implies the one for $M$.

(ii)  Define the index $(M:N)$ of any two torsion-free $\mathcal O_F$-modules $M$ and $N$ of rank $d$ inside $V$ by $(M:N) = (M:P)/(N:P)$ for any submodule $P$ 
of finite index contained in $M$ and $N$.  Let $\{x_1, \cdots, x_d\}$ be a basis of $V= M\otimes_{\mathcal O_F}F$.  
Let $\{\epsilon_j\}$ be the standard basis of $ \mathbb C^{r_1+2r_2}$. Let $\{y_1, \cdots, y_d\}$ be a $F$-basis of $M^{\vee}\otimes_{\mathcal O_F}F$ 
dual to $\{x_1, \cdots, x_d\}$. Then $(<x_i\otimes\epsilon_j>, (M: <x_1, \cdots, x_d>))$ is a representative for 
the $\mathcal O_F$-integral structure $j_M$.  Its $\mathbb Z$-dual is $(<y\otimes\epsilon_j>, (<x_1, \cdots, x_d> : M))$. But $(<y\otimes\epsilon_j>, (M^{\vee}: <y\otimes\epsilon_j>))$ is a representative for $j_{M^{\vee}}$ and since $(M:N)^{-1} = (M^{\vee}, N^{\vee})$, we are done. 
\end{proof} 

\subsubsection{Pairings and determinants}\label{pairings} Let $N$ be a finitely generated abelian group and $\psi: N/{N_{tor}} \times N/{N_{tor}} \to \mathbb R$ be a non-degenerate symmetric bilinear form on $N$.  One defines $\Delta(N) \in \mathbb C^{\times}/{\{\pm1\}}$ as  the determinant of the matrix $\psi(b_i, b_j)$ divided by  $(N: N_0)^2$ where $N_0$ is the subgroup of finite index generated by a maximal linearly independent subset $b_i$ of $N$; so $\Delta(N)$ is independent of the choice of $b_i$'s and incorporates the order of the torsion subgroup of $N$. We can rewrite $\psi$ as the integral structure $(C_{\bullet}, L_{\bullet})$ with $L_0=N, L_1= \textrm{Hom}(N, \mathbb Z)$: 
\[C_0= N\otimes\mathbb C \xrightarrow{\psi_*}  \textrm{Hom}(N, \mathbb Z)\otimes\mathbb C= C_1, \qquad  \Delta(N) = \frac{det(\psi_*) }{[N_{tor}]^2}; \]
here $det(C_{\bullet}, L_{\bullet})$ is the determinant $\det(\psi_*)$ of $\psi$ computed with the bases $L_0$ and $L_1$. 
Given a short exact sequence of finitely generated groups $ 0 \to N' \to N \to N'' \to 0$
which splits over $\mathbb Q$ as an orthogonal direct sum 
\[N_{\mathbb Q} \cong N'_{\mathbb Q} \oplus N''_{\mathbb Q}\]
 with respect to a definite pairing $\psi$ on $N$, one has the following standard relation
\begin{equation}\label{deltann} \Delta(N) = \Delta(N')~\cdot~\Delta(N'').
\end{equation}

\subsection{Motives} 
 We start with the integral structures underlying the Betti and de Rham realizations of the motives $M= h^j(V_0)(r)$ for $0\le j \le 2(d-1)$.

\begin{definition}\label{betti}  The integral structure $H^j_B(V_{\mathbb C}, \mathbb Z(r))^+$ is defined as the complex vector space $H^j_B(V_{\mathbb C}, \mathbb C(r))^{+}$  together with the homomorphism
\[H^j_B(V_{\mathbb C}, \mathbb Z)^+ \to H^j_B(V_{\mathbb C}, \mathbb C)^{+}\xrightarrow{\times (2\pi i)^r} H^j_B(V_{\mathbb C}, \mathbb C(r))^{+} ;\]
the first map is the natural map and the second map is multiplication by $(2\pi i)^r$.\end{definition}
 
 \begin{example} When $V=S$, one has $H^0(S_{\mathbb C}, \mathbb Z(r))^{+}$ is isomorphic to $\mathbb Z^{r_2}$ when  $r$ is odd and $\mathbb Z^{r_1+r_2}$ when $r$ is even.
\end{example} 

 \subsection{de Rham realization of $M= h^j(V_0)(r)$} The de Rham realization $M_{dR} = H^j_{dR}(V_{\mathbb C}, \mathbb C(r))$ is $H^j_{dR}(V_{\mathbb C}, \mathbb C)$ with the Hodge filtration shifted by $r$, i.e., the Tate twist in de Rham realization shifts the Hodge filtration but does not change the vector space. The tangent space $(t_M)_{\mathbb C}$ of $M$ is defined as follows
 \[(t_M)_{\mathbb C} = \frac{H^j_{dR}(V_{\mathbb C}, \mathbb C(r))}{F^0} = \prod_{k<r} H^{j-k}(V_{\mathbb C}, \Omega^k_{V_{\mathbb C}}) = \prod_{\sigma \in \tilde{S} }~ \prod_{k<r} H^{j-k}(V_{\sigma}, \Omega^k_{V_{\sigma}}).\] 
 
 Let $\lambda^k\Omega_{V/S}$ denote the $k$th derived exterior power of the sheaf $\Omega_{V/S}$ on $V$; see the appendix of \cite{SL2020}.  
 \begin{definition}\label{integral-dr-defn} The integral structure on the complex vector space $(t_M)_{\mathbb C}$ is given by the map
 \[t_M:= \prod_{k<r} H^{j-k}(V, \lambda^k \Omega_{V/S})  \to \prod_{\sigma \in \tilde{S}}~ \prod_{k<r} H^{j-k}(V_{\sigma}, \Omega^k_{V_{\sigma}}) = (t_M)_{\mathbb C}.\] 
 \end{definition}
As in (\ref{im-integral}), one views $H^{j-k}(V, \lambda^k \Omega_{V/S}) $ as an abelian group; its structure as a module over $\mathcal O_F$ is not used.

\begin{remark} 
We can use the Hodge decomposition of $(t_M)_{\mathbb C}$ to decompose the period map $\alpha_M$ of (\ref{alphaM})  as a sum of $\alpha_k$ where $\alpha_k$ is the composite map 
\[H^j(V_{\mathbb C}, \mathbb C(r))^{+} \xrightarrow{\alpha}  (t_M)_{\mathbb C}  = \prod_{k<r} H^{j-k}(V_{\mathbb C}, \Omega^k_{V_{\mathbb C}}) \to H^{j-k}(V_{\mathbb C}, \Omega^k_{V_{\mathbb C}}).\] The enhanced period map $\gamma_M$ of \cite[\S 2.3]{SL2020} 
\[\gamma_M: H^j(V_{\mathbb C}, \mathbb C(r))^{+}  \to (t_M)_{\mathbb C},\] likewise, decomposes as a sum of  $\gamma_k:= \Gamma^*(r-k)~\cdot~\alpha_k$. Here $\Gamma$ is the classical gamma function and $\Gamma^*(m)$ denotes the special value at $s=m$. \qed
 \end{remark}

\subsection{Integral structures on certain motivic exact sequences}\label{conj-integral}
The conjectural formula (\ref{conj}) from  \cite[Conjecture 3.1]{SL2020}  predicts that the special value $\zeta^*(V,r)$ of $\zeta(V,s)$ at $s=r$ is equal to  $ \pm \chi(V,r)$ (up to powers of two) defined as follows
\[ \chi(V,r)= \frac{\chi_{A,C}(V,r)}{\chi_B(V,r)} = \frac{\chi(C(r))~\cdot~\chi_A(V,r)~\cdot~\chi_{A'}(V,r)}{\chi_B(V,r)}.\]
In more detail, 
\begin{itemize}
\item  $M_j$ is the motive $h^j(V_0,  \mathbb Z(r))$ with $0\le j \le 2d-2$ with the enhanced period map
\[\gamma_{M_j} : H^j_B(V_{\mathbb C}, \mathbb Z(r))^{+}_{\mathbb C} \to t_{M_j}\]
\item we fix an integral structure on $ \textrm{Ker}(\gamma_{M_j})$ and $ \textrm{Coker}(\gamma_{M_j})$ for all $j$ and use them implicitly in all that follows. Under this condition, the final conjecture is actually independent of the choice of these integral structures. 
\item the exact sequences $B(j,r)$ with integral structures 
 \begin{equation}\label{Bjr} 
 0 \to \textrm{Ker}(\gamma_{M_j}) \to H^j_B(V_{\mathbb C}, \mathbb Z(r))^{+}_{\mathbb C} \to t_{M_j} \to \textrm{Coker}(\gamma_{M_j}) \to 0.
\end{equation} 
We write 
\[\chi_B(V,r) =\prod_j (\chi(B(j,r)))^{(-1)^j} = \frac{\chi(B(0,r))~\cdot~\chi(B(2,r))~\cdots}{\chi(B(1,r))~\cdot~\chi(B(3,r)) \cdots}.\]
\item  the integral structures $A(j,r)$ for $0\le j\le \textrm{min}(2d-1, 2r-3)$ and $A'(j,r)$ for $2d+1\ge j\ge \textrm{max}(0, 2r+1)$ defined as
\begin{equation}\label{Ajr}
A(j,r): H^{j+1}_W(V, \mathbb Z(r))_{\mathbb C} \to \textrm{Coker}(\gamma_{M_j}), \qquad A'(j,r): \textrm{Ker}(\gamma_{M_j}) \to H^{j+2}_W(V,\mathbb Z(r))_{\mathbb C}.
\end{equation} 
We write 
\begin{equation}\label{aaprime}
\chi_A(V,r) = \prod_{j=0}^{\textrm{min}(2d-1, 2r-3)} (\chi(A(j,r)))^{(-1)^j}, \qquad \chi_{A'}(V,r) = \prod_{j=\textrm{max}(0, 2r+1) }^{j=2d-1} (\chi(A'(j,r)))^{(-1)^j}.
\end{equation}
The map in $A(j,r)$ is Beilinson's regulator from (\ref{FPles}) , usually stated as a map from algebraic K-theory of $V_0$ to the Deligne cohomology of $V_0$; it is conjectured to be an isomorphism. The group $\textrm{Ker}(M_{2d-2-j, d-r})$ can be identified with the dual of $\textrm{Coker}(M_{j,r})$ - see \cite[\S 2]{SL2020} or \cite[\S 6.9]{MR1206069}. The integral structures $A'(j,r)$ are defined as the dual of $A(2d-2-j,d-r)$ in \cite{SL2020}. 
\item the exact sequence $C(r)_{\mathbb C}$ with integral structure (for $0 \le r\le d$)
\begin{multline}\label{Cr} 
0 \to H^{2r-1}_W(V, \mathbb Z(r))_{\mathbb C} \to \textrm{Coker}(\gamma_{M_{2r-2}}) \to H^{2r+1}_W(V, \mathbb Z(r))_{\mathbb C} \xrightarrow{h}\\ 
\to H^{2r}_W(V, \mathbb Z(r))_{\mathbb C} \to \textrm{Ker}(\gamma_{M_{2r}}) \to H^{2r+2}_W(V, \mathbb Z(r))_{\mathbb C} \to 0.\end{multline} 
The maps here are integral versions of the maps from  (\ref{FPles}): the first is Beilinson's regulator, the second is the dual of the cycle map, the third is the Arakelov intersection/height pairing map, the fourth map is the  cycle class map, and the fifth is the dual of Beilinson's regulator. Beilinson's conjectures \cite[Conjectures 2.3.4-2.3.7]{SL2020} imply that $C(r)_{\mathbb C}$ is exact; we refer to \cite[\S 2]{SL2020} for more details.  Integral versions of (\ref{FPles}) have also been introduced by Flach-Morin  \cite[Conjecture 2.9 and (28), (29)] {MR3874942}.
\end{itemize} 

\begin{table}[h!]
  \begin{center}
    \caption{Tabulation of invariants in the various integral structures for the motives $M_j =h^j(V_0)(r)$}
    \label{tab:table2}
    \begin{tabular}{l|c|c} 
      \textbf{Groups} & {\bf Type $B$ }& {\bf Type $A$ or $C$}\\
      \hline
            $H^j_W(V, \mathbb Z(r))$ and $\textrm{max}(1, 2r+2) < j < 2d+2$ & - & $A'(j-2,r)$ \\
      \hline
      $H^j_W(V, \mathbb Z(r))$ and $0< j< \textrm{min}(2d+1,2r-2)$ & - & $A(j-1,r)$ \\
      \hline
      $H^j_W(V, \mathbb Z(r))$ and $2r-2 < j < 2r+3$ & - & $C(r)$ \\
      \hline
      $\textrm{Ker}(\gamma_{M_j})$ and $j \ge 2r+1$ &  $B(j,r)$ & $A'(j,r)$\\
      \hline
      $\textrm{Ker}(\gamma_{M_j}) =0$ for $j \le 2r-1$ &  $B(j,r)$ & -\\
      \hline
      $\textrm{Coker}(\gamma_{M_j}) =0$ for $j > 2r$& $B(j,r)$  & - \\
      \hline
       $\textrm{Coker}(\gamma_{M_j}) =0$ for $j = 2r$& $B(j,r)$  & - \\
      \hline
       $\textrm{Coker}(\gamma_{M_j})$ and $j\le 2r-3$& $B(j,r)$  & $A(j,r)$ \\
      \hline
        $\textrm{Ker}(\gamma_{M_{2r}})$ & $B(2r,r)$  & $C(r)$ \\
      \hline
      $\textrm{Coker}(\gamma_{M_{2r-2}})$& $B(2r-2,r)$  & $C(r)$\\
      \hline
      $\textrm{Coker}(\gamma_{M_{2r-1}}) =0 =\textrm{Ker}(\gamma_{M_{2r-1}})  $ &  &\\
      \hline
    \end{tabular}\label{table}
  \end{center}
\end{table}

\begin{remark} (i) By \cite[\S 7.1]{MR1206069} or \cite[\S 2.1 and (2.2.1)]{MR1265544}, if the weight of $M$ is positive (so $j>2r$), then one has $\textrm{Coker}(\gamma_M) =0$. Dually, if $j-2r <0$, then one has $\textrm{Ker}(\gamma_M) =0$. It remains to consider the case $M= h^{2r}(V_0, \mathbb Z(r))$. Since $N = h^{2r}(V_0, \mathbb Z(r+1))$ has negative weight, it follows that $\textrm{Ker}( \gamma_N)=0$.  Hence  $\textrm{Coker}(\gamma_{M_{2r}}) =0$ for as it is dual to  $\textrm{Ker}( \gamma_N)$.

(ii) Note that each group $H^i_W(V,\mathbb Z(r))$ appears exactly once.  

All the exact sequences above are variants of motivic exact sequences (\ref{FPles})  but equipped with an integral structure.  
 Our convention is that the terms $ \textrm{Ker}(\gamma_{M})$ has even degree and  $\textrm{Coker}(\gamma_{M})$ has odd degree.  If $j=2r-1$, then $\textrm{Ker}(\gamma_{M_j})=0 = \textrm{Coker}(\gamma_{M_j})$; if $j \neq 2r-1$, these terms occur exactly twice in the conjecture, once in the $B$-complexes and once in $A$, $A'$ and $C$-complexes. Since $\chi_B$ is the denominator and the others ($\chi_A$, $\chi_{A'}$ and $\chi(C)$) are in the numerator of $\chi(V,r)$, the conjecture is independent of the choice of integral structures on the terms $ \textrm{Ker}(\gamma_{M})$ and  $\textrm{Coker}(\gamma_{M})$. 
 
 (iii) The conjectural formula (\ref{conj}) from  \cite[Conjecture 3.1]{SL2020}  tacitly assumes the
  generalized Beilinson-Soul\'e conjecture from \cite{SL2020} which predicts 
 \[ \textrm{if}~r<0,~\textrm{then}~H^i_W(X, \mathbb Z(r)) = \textrm{finite~2-group}~\textrm{for}~i\le 1.\] If this is not assumed, then (\ref{aaprime}) has to be replaced with 
\[\chi_A(V,r) = \prod_{j=0}^{2r-3} (\chi(A(j,r)))^{(-1)^j}, \qquad \chi_{A'}(V,r) = \prod_{j=2r+1}^{j=2d-1} (\chi(A'(j,r)))^{(-1)^j}\]
and the limits on $j$ in the first two lines of Table \ref{tab:table2} have to be replaced with 
$2d+2> j>2r+2$ (in the first line) and  $0< j< 2r-1$ (in the second line). We refer to \cite{SL2020} for more details.

 (iv) We refer to \cite{SL2020} for the conjectures such as the finite generation of $H^*_W(V, \mathbb Z(r))$ which are implicit in the formulation of (\ref{conj}). \qed
 \end{remark} 
 It is instructive to consider the conjecture \cite[Conjecture 3.1]{SL2020} in the case of number fields before delving into the case of arithmetic surfaces.  
\subsection{The case of number fields: $V=S$, $d=1$, and $r=0$ and $r=1$}\label{d=1} This example is worked out in \cite[\S7]{SL2020}.  

 \subsubsection{The case $r=0$} The motives are $M_j =h^j(S)(0)$. 
 
 We know $H^j_{et}(S, \mathbb Z(1)) =0$ if $j <1$, that $H^1_{et}(S, \mathbb Z(1)) = \mathcal O_F^{\times}$, $H^2_{et}(S, \mathbb Z(1)) = \textrm{Pic}(S)$, $H^3_{et}(S, \mathbb Z(1)) =0$ (up to a finite 2-group), and $H^0_{et}(S, \mathbb Z(0)) = \mathbb Z$. 
 
 It immediately follows from the definitions that $H^0_W(S, \mathbb Z(0)) =\mathbb Z$, $H^1_W(S, \mathbb Z(0)) =0$ (up to a finite $2$-group), and 
 $H^3_W(S, \mathbb Z(0)) = \mu_F^{\vee}$, the dual of the roots of unity in $F$.
 
 There is also an exact sequence
 \[0 \to \textrm{Pic}(S)^{\vee} \to  H^2_W(S, \mathbb Z(0)) \to \textrm{Hom}(\mathcal O_F^{\times}, \mathbb Z) \to 0.\]
 
Note that  $t_{M_j}=0$ for all $j$. 

$\chi_B(S,0)=1$: The sequence $B(j,0)$ is zero for all $j \neq 0$ and we can arrange $\chi(B(0,0)) =1$ by using the isomorphism 
\[ B(0,0):\qquad \textrm{Ker} (\gamma_{M_0}) \xrightarrow{\sim} H^0_B(S_{\mathbb C}, \mathbb Z(0))_{\mathbb C}^+\]
to transfer the integral structure $H^0_B(S_{\mathbb C}, \mathbb Z(0))^+$ to   
$\textrm{Ker}(\gamma_{M_0})$.

Since $A(j,r)$ contributes only for $j \le 2r-3 =-3$, it follows that there is no contribution from any of the $A(j,0)$-terms.  
 We see that $A'(j,0) =0$ unless $j=1$, in which case $A'(1,0) = \mu_F^{\vee}$ in degree two. Thus, we obtain $\chi(A'(1,0)) = w$ and 
 \[\chi_{A'}(S,0) = \frac{1}{w}.\] 
  The sequence $C(0)_{\mathbb C}$  for $V=S$
 \[0\to H^0_W(S, \mathbb Z(0))_{\mathbb C} \to \textrm{ker}(\gamma_{M_0}) \to H^2_W(S, \mathbb Z(0))_{\mathbb C} \to 0\] becomes the sequence (our convention is that $\textrm{Ker}(\gamma)$ has even degree, say, degree two)
 \begin{equation}\label{conj-S-zero} 0 \to \underset{\textrm{degree one}}{\mathbb C} \to \underset{\textrm{degree two}}{\mathbb C^{r_1+r_2}} \to \underset{\textrm{degree three}}{\textrm{Hom}(\mathcal O_F^{\times}, \mathbb C)} \to 0, 
  \end{equation} 
 with $H^2_W(S, \mathbb Z(0))$ providing the standard integral structure $\textrm{Hom}(\mathcal O_F^{\times}, \mathbb Z)$ on the last term and the second map being the dual of the classical regulator. So $\chi(C(0)) = hR$, and so 
 \[ \chi(S,0) =  \frac{\chi_{A,C}(S,0)}{\chi_B(S,0)} = \frac{\chi(C(0))~\cdot~\chi_{A'}(S,0)}{1} = \frac{hR}{w}.\] As is well known, 
 \begin{equation}\label{cnf0} \zeta^*(S,0) =   -hR/{w}.\end{equation}  So the conjectural formula (\ref{conj}) from  \cite[Conjecture 3.1]{SL2020} is valid in this case.
 \subsubsection{The case $r=1$} The motives are $M_j=h^j(S)(1)$. 
 
 As $A(j,1)$ contribute only for $j <-1$, we have $\chi_A(S,1) =1$.   
 The complex $C(1)_{\mathbb C}$ (our convention is that $\textrm{Coker}(\gamma)$ has odd degree, say one)
 \[0 \to H^{1}_W(S, \mathbb Z(1))_{\mathbb C} \to \textrm{Coker}(\gamma_{M_{0}}) \to H^{3}_W(S, \mathbb Z(1))_{\mathbb C} \to 0 \to 0 \cdots \]
 becomes
 \begin{equation}\label{conj-S-1}0 \to \underset{\textrm{degree zero }}{\mathcal O_F^{\times}\otimes\mathbb C} \to \underset{\textrm{degree one}}{\mathbb C^{r_1+r_2}} \to \underset{\textrm{degree two }}{\mathbb C} \to 0,\end{equation}
 with the last two terms getting the standard basis and the first term a basis from $\mathcal O_F^{\times}$. This is exactly as in \cite[\S 7]{SL2020}. Thus, $det(C(1))$ is the classical regulator $R$. One checks that $H^i_W(S, \mathbb Z(1)) =0$ for $i> 3$ and that $H^3_W(S, \mathbb Z(1)) = \textrm{Hom}(H^0_{et}(S, \mathbb Z(0)),  \mathbb Z)$.  As $H^i_W(S,\mathbb Z(1)) = H^i_{et}(S, \mathbb Z(1))$ for $i \le 2$, this completes the computation of  $H^*_W(S, \mathbb Z(1))$. From this, we see that
 \[\chi (C(1)) = \frac{hR}{w}.\]
 
 We need to consider $A'(j,1)$ only when $j \ge 3$. Since $t_{M_j} =0$ for $j>0$ and $H^{j+2}_W(S, \mathbb Z(1)) =0$ for $j\ge 3$, we see $A'(j,1) =0$ for $j\ge 3$.  As a result, $\chi_{A'}(S,1)=1$.
 
 Finally, $B(j,1) =0$ if $j \neq 0$ and $B(0,1)$ given by
 \[0\to  H^0_B(S_{\mathbb C}, \mathbb Z(1))^{+}_{\mathbb C} \to t_{M_0} \to \textrm{Coker}(\gamma_{M_0}) \to 0\] becomes (our convention is that $\textrm{Coker}(\gamma)$ has odd degree)
 \[0 \to \underset{\textrm{degree one }}{\mathbb C^{r_2}} \to \underset{\textrm{degree two }}{\mathcal O_F\otimes_{\mathbb Z}\mathbb C} \to \underset{\textrm{degree three }}{\mathbb C^{r_1+r_2}} \to 0.\]
 The first map is the composition of the natural inclusion of $\mathbb C^{r_2}$ in $\mathbb C^{r_1+2r_2}$ (multiplied by $2\pi i$) with the inverse of $\psi$ of (\ref{standard-psi}).  The integral structure on $\textrm{Coker}(\gamma_{M_0}) \cong \mathbb C^{r_1+r_2}$ is defined by $H^0_B(S, \mathbb Z(1))^{-} \cong \mathbb Z^{r_1+r_2}$. Since the determinant of $\psi$ with respect to the usual bases is $\sqrt{d}_F$ and the integral structures are torsion-free, we obtain
\[ det (B(0,1)) = \frac{\sqrt{d_F}}{(2\pi i)^{r_2}} = \chi(B(0,1)).\]

Hence we obtain
\[\chi(S,1) =\frac{\chi(C(1))}{\chi(B(0,1))} = \frac{hR(2\pi i)^{r_2}}{w\sqrt{d_F}};\]
as the sign of $d_F$ is $(-1)^{r_2}$, this is equal to the usual formula (up to a power of two)
\begin{equation}\label{cnf1} \zeta^*(S,1) = \frac{2^{r_1}h R(2\pi)^{r_2}}{w\sqrt{|d_F|}}.\end{equation}
So the conjectural formula (\ref{conj}) from  \cite[Conjecture 3.1]{SL2020} is valid in this case too.

\subsection{The case of arithmetic surfaces: $V=X$ and $r=1$} We now turn to an explicit description of the terms that enter into the description of $\zeta^*(X,1)$; this uses the motives $M_j = h^j(X_0)(1)$ for $j=0,1,2$ of the algebraic curve $X_0$ over $F$ and $t_{M_j} = H^j(X_0, \mathcal O_X)$.  We begin with the Weil-\'etale motivic cohomology groups of $X$. 

\subsubsection{The groups $H^*_W(X,\mathbb Z(1))$}\label{Euler1} Fix an arithmetic surface $\pi:X\to S$.
 
Since the motivic complex $\mathbb Z(1)$ is $ \mathbb G_m[-1]$, one has the identifications
\[H^j_{et}(X, \mathbb Z(1)) = H^{j-1}_{et}(X, \mathbb G_m),\]  
 the Picard group $\textrm{Pic}(X) = H^2_{et}(X, \mathbb Z(1))$, and the Brauer group $\textrm{Br}(X) = H^3_{et}(X, \mathbb Z(1))$.  

Recall that  $H^i_W(X, \mathbb Z(1))$ \cite[\S 2.1]{SL2020} are defined as \'etale motivic cohomology $H^i_{et}(X, \mathbb Z(1))$ for $i \le 2$ and, for $i >2$, as the dual of \'etale motivic cohomology  
\[H^i(\textrm{RHom}(R\Gamma_{et}(X, \mathbb Z(1)), \mathbb Z [-3])),\] thereby sitting in an exact sequence (for $i >2$) 
\[0 \to \textrm{Ext}^1(H^{6-i}_{et}(X, \mathbb Z(1)), \mathbb Z) \to H^i_W(X, \mathbb Z(1)) \to \textrm{Hom}(H^{5-i}_{et}(X, \mathbb Z(1)), \mathbb Z) \to 0.\]

 As $H^i_W(X, \mathbb Z(1)) =0$ for $i=0$ and $i \ge 6$, the following is a complete description of $H^*_W(X, \mathbb Z(1))$: 
\begin{equation}\label{weil2} 
 H^1_W(X, \mathbb Z(1)) = {{\mathcal O}_F}^{\times}, \quad H^2_W(X, \mathbb Z(1)) = {\rm Pic}(X), \quad H^5_W(X, \mathbb Z(1)) = {\rm Ext}^1({\mathcal O_F}^{\times}, {\mathbb Z}),
\end{equation} 

 \begin{equation}\label{weil3} 0 \to {\rm Ext}^1({\rm Br}(X), {\mathbb Z}) \to H^3_W(X, \mathbb Z(1)) \to {\rm Hom}({\rm Pic}(X), \mathbb Z) \to 0,
 \end{equation} 

\begin{equation}\label{weil4}  0 \to {\rm Ext}^1({\rm Pic}(X), {\mathbb Z}) \to H^4_W(X, \mathbb Z(1)) \to {\rm Hom}(\mathcal O_F^{\times}, \mathbb Z) \to 0.\end{equation}

\begin{remark*} As the finite generation of $\textrm{Pic}(X)$ is well known (theorem of Mordell-Weil-Roquette), the  finite-generation of $H^*_W(X,\mathbb Z(1))$ reduces to the finiteness of $\textrm{Br}(X)$.\qed \end{remark*}  

 \subsubsection{The complexes $A(j,1)$ and $A'(j,1)$}\label{aj} In our case, only $A'(3,1)$ could be nonzero.

In our case, $r=1$. For $A(j,1)$ to intervene, one needs $j\le 2r-3 =-1$ as $2d-1=3$. Since $H^0_W(X,\mathbb Z(1)) =0$, we obtain $\chi(A(j,1)=1)$ for all $j$.   

For $A'(j,1)$ to intervene, one needs $j \ge 3$. For $j=3$, since  $H^5_W(X,\mathbb Z(1))$ is a finite group of order $w$, we have  $det(A'(3,1)) =1$ and  $\chi(A'(3,1)) =w$. If $j>3$, then $H^{j+2}_W(X, \mathbb Z(1)) =0$ and $M_j =0$; so $\chi(A'(j,1)) =1$ for $j>3$. Thus,
\begin{equation}\label{chiaa'} \chi_A(X,1) =1, \qquad \chi_{A'}(X,1) = \frac{1}{w}.\end{equation} 

\subsubsection{The groups $H^i(X, \mathcal O_X)$}
 Since $\pi:X \to S$ is proper, $H^i(X, \mathcal O_X)$ are finitely generated $\mathcal O_F$-modules. 
\begin{lemma}\label{h2oh}   $H^i(X, \mathcal O_X)$ is zero for $i\ge 2$.
\end{lemma}
\begin{proof}  Since $X\to S$ is a relative curve, $H^i(X_s, \mathcal O_{X_s}) =0$ (for $i >1$) at all points $s\in S$. Since $S$ is a reduced noetherian scheme and $\pi: X\to S$ is proper and $\mathcal O_X$ is a coherent sheaf on $X$ flat over $S$, by Grothendieck's theorem on formal functions (see L.~Illusie's article \cite[8.5.18]{FGA}), $R^i\pi_*\mathcal O_X =0$ for $i>1$.\end{proof}

Using \S \ref{aj} and Lemma \ref{h2oh}, we find that 
\begin{equation}\label{chix1} \chi(X,1) = \frac{{\chi (C(1))}}{\chi(A'(3,1))}~\cdot~\frac{\chi(B(1,1))}{\chi(B(0,1))} =\frac{{\chi (C(1))}}{w}~\cdot~\frac{\chi(B(1,1))}{\chi(B(0,1))},\end{equation} 
 where  the exact sequence $B(j,1)$ is
 \begin{equation}\label{Bj1} 
 0 \to \textrm{Ker}(\gamma_{M_j}) \to H^j_B(X_{\mathbb C}, \mathbb Z(1))^{+}_{\mathbb C} \to t_{M_j} \to \textrm{Coker}(\gamma_{M_j}) \to 0.
\end{equation} 
 and $C(1)_{\mathbb C}$ is
\begin{equation}\label{C1} 
0 \to H^1_W(X, \mathbb Z(1))_{\mathbb C} \to \textrm{Coker}(\gamma_{M_0}) \to H^3_W(X, \mathbb Z(1))_{\mathbb C} \to H^2_W(X, \mathbb Z(1))_{\mathbb C} \to \textrm{Ker}(\gamma_{M_2}) \to H^4_W(X, \mathbb Z(1))_{\mathbb C} \to 0.\end{equation} 
So Conjecture \ref{main conjecture} provides an intrinsic formula for $\zeta^*(X,1)$ using only  $B(0,1)$, $B(1,1)$, $C(1)$, and $H^*_W(X, \mathbb Z(1))$. We compute $\chi(B(0,1))$ and explain how $\chi(B(2,1))$ can be assumed to be one.  

\subsubsection{$\chi(B(0,1))$} This has been computed (\cite[\S 7]{SL2020}) in \S \ref{d=1}. Here
\begin{equation}\label{gamma-m0}  
\textrm{Ker}(\gamma_{M_0}) =0, \qquad \textrm{Coker}(\gamma_{M_0}) \cong \mathbb C^{r_1+r_2} ,\qquad H^0(X, \mathcal O) = \mathcal O_F.
\end{equation}  Our convention for $B(0,1)$ is that $\textrm{Coker}(\gamma_M)$ is in odd degree, say degree three.
As the integral structure
\[B(0,1):\qquad  0 \to  \underset{\textrm{degree one }}{H^0_B(X_{\mathbb C}, \mathbb Z(1))^{+}_{\mathbb C}} \to  \underset{\textrm{degree two}}{(\mathcal O_F)_{\mathbb C}} \to  \underset{\textrm{degree three}}{\textrm{Coker}(\gamma_{M_0})} \to 0\]
is torsion-free, we have (again using that the sign of $d_F$ is $(-1)^{r_2}$)
\begin{equation}\label{bzero1} \chi (B(0,1)) = det (B(0,1)) = \frac{\sqrt{d_F}}{(2\pi i)^{r_2}} =  \frac{\sqrt{|d_F|}}{(2\pi)^{r_2}}.\end{equation}
\subsubsection{The sequence $B(2,1)$} By Lemma \ref{h2oh},  $t_{M_2} = H^2(X_0, \mathcal O) =0$; so, in $B(2,1)$
\[ \textrm{Ker}(\gamma_{M_2}) \xrightarrow{\sim} H^2_B(X_{\mathbb C}, \mathbb Z(1))^{+}_{\mathbb C},\] we use the isomorphism to transfer the integral structure from $H^2_B(X_{\mathbb C}, \mathbb Z(1))^{+}$ to $\textrm{Ker}(\gamma_{M_2})$, thereby obtaining 
\begin{equation}\label{btwo1}\chi(B(2,1)) = 1.\end{equation}

\section{The conjecture of Birch-Swinnerton-Dyer-Tate}\label{Background} 

\subsection{Statement of the conjecture \cite{Tate1966, Gross82}} Our presentation of this conjecture follows (almost verbatim) \cite[\S~2]{Gross82} where it is formulated in terms of N\'eron models.

  Let $A$ be an abelian variety of dimension $g$ over ${\rm Spec}~F$.  Let $L(A,s)$ be its L-series over $F$, which is defined -- see (\ref{LAV}) -- by an Euler product convergent in the half-plane Re~$(s)> \frac{3}{2}$. We assume that $L(A,s)$ has an analytic continuation to the entire complex plane. we write $L^*(A,1)$ for its special value at $s=1$. Writing
\[ L(A,s)~\sim~~c~\cdot~(s-1)^{r(A)} \qquad s\to 1,\]
where $r(A)$ is a non-negative integer and $c$ is non-zero (the special value $L^*(A,1)$). Birch and Swinnerton-Dyer and Tate have conjectured that $r(A)$ is equal to the rank of the finitely generated group $A(F)$ of points of $A$ over $F$; they have also given a conjectural formula for $L^*(A,1)$ in terms of certain arithmetic invariants of $A$ which we now recall.

Let $\mathcal N$ be the N\'eron model of $A$ over $S$. Let $\mathcal N^0$ be the largest open sub-group scheme of $\mathcal N$ in which all fibers are connected.  For any non-zero prime $v$ of $\mathcal O_F$ (i.e., $v \in S$), we write 
\[\mathcal N_v, \qquad \mathcal N^0_v\] for the fibers over $v$. Thus, $\mathcal N_v$ is a commutative group scheme over $k(v)$ and $\mathcal N^0_v$ is the connected component of the identity of $\mathcal N_v$. We have an exact sequence of group schemes over $k(v)$:
\[0 \to H_v \to \mathcal N_v^0 \to B_v \to 0\]
where $H_v$ is a commutative linear group scheme and $B_v$ is an abelian variety. For almost all $v$, we have 
\[\mathcal N_v =\mathcal N^0_v = B_v;\]these are the primes where $A$ has good reduction. 
Let $Y_v = \textrm{Hom}(H_v, \mathbb G_m)$ be the character group of $H_v$ and $\pi_v$ be the Frobenius endomorphism of $B_v$. 

We write $\Phi_v=\pi_0(\mathcal N_v)$ for the finite group scheme of connected components of the fiber $\mathcal N_v$ over $v$. The Tamagawa number of $A$ at $v$, denoted $c_v$, is the size of the group $\Phi_v(k(v))$ of the group of rational points of $\Phi_v$. The Tamagawa number $c_v$ is one for almost all $v \in S$; we may define the product 
\begin{equation}\label{tamagawas} P_{A, fin} = \prod_{v\in S} c_v.\end{equation}

Let $\Gamma_v$ be a decomposition group for $v$ in $\Gamma_F = \textrm{Gal}(\bar{F}/F)$ and let $I_v$ be the inertia subgroup
of $\Gamma_v$ and $\sigma_v$ be an arithmetic Frobenius, a topological generator of the quotient $\Gamma_v/{I_v}$.  For any prime $\ell$ distinct from the characteristic of $k(v)$, consider the $\ell$-adic Tate module $T_{\ell}A$ of $A$; this is a free $\mathbb Z_{\ell}$-module of rank $2g$ which admits a continuous $\mathbb Z_{\ell}$-linear action of $\Gamma_F$. We define the local $L$-factor of $A$ at $v$ by the formula
\begin{equation}\label{localfactor} L_v(A,t) = \textrm{det}(1-\sigma_v^{-1}t~|~\textrm{Hom}_{\mathbb Z_{\ell}}( T_{\ell}A, \mathbb Z_{\ell})^{I_v}).\end{equation}
The characteristic polynomial $L_v(A,t)$ has integral coefficients which are independent of $\ell$: this is evident from the formula
\[L_v(A,t) = \textrm{det}(1-\sigma_v t~|~ Y_v)~\cdot~\textrm{det}(1-\pi_v t~|~ T_{\ell} B_v )\]
where the first determinant is clearly independent of $\ell$ and the second is the characteristic polynomial of an endomorphism of the abelian variety $B_v$ which has integral coefficients independent of $\ell$ by a theorem of A.~Weil.  

The roots of the first factor have complex absolute value one whereas those of the second factor have complex absolute value $q_v^{-1/2}$. Therefore, the global $L$-series $L(A,s)$ defined by
\begin{equation}\label{LAV} L(A,s) = L(A/F,s) = \prod_{v\in S} \frac{1}{L_v(A, q_v^{-s})}\end{equation}
converges for Re~$(s) >\frac{3}{2}$. We shall assume that it has an analytic continuation to the entire complex plane. 

\subsection{The archimedean period of $A$}\label{arch-periodes}  Let $\omega_{\mathcal N}$ denote the {\it projective} $\mathcal O_F$-module of invariant differentials of the N\'eron model $\mathcal N$. Since $g$ is the rank of $\omega_{\mathcal N}$, the module $\Lambda^g\omega_{\mathcal N}$ is a rank one $\mathcal O_F$-submodule of $H^0(A, \Omega^g)$ - a vector space over $F$ of dimension one. Let $\{w_1, \cdots, w_g\}$ be a $F$-basis of $H^0(A, \Omega^1)$ and put $\eta = w_1\wedge\cdots \wedge w_g$. We have
\begin{equation}
\Lambda^g\omega_{\mathcal N} = \eta~\cdot~\mathfrak a_{\eta} \quad \subset H^0(A, \Omega^g) 
\end{equation} 
where $\mathfrak a_{\eta}$ is a fractional ideal of $F$ and 
\begin{equation}\label{c-eta} \mathfrak a_{\eta} = \mathfrak a_{c\eta} ~\cdot~(c) \qquad c \in F^{\times}.\end{equation}

For any $\sigma \in S_{\mathbb R}$ (corresponds to a real embedding $\sigma: F \to \mathbb R$), let $H^{+}$ denote the submodule of $H_1(A_{\sigma}(\mathbb C), \mathbb Z)$ which is fixed by complex conjugation; it is a free $\mathbb Z$-module of rank $g$ and if $\{ \gamma_1, \cdots, \gamma_g\}$ denotes a basis, we put
\begin{equation}\label{grossreal}
p_{\sigma} (A, \eta) = [(\pi_0(A_{\sigma}(\mathbb R)))]~\cdot~\Bigg|~\textrm{det} \bigg( \bigg(\int_{\gamma_i}\omega_j \bigg)\bigg)~\Bigg|.\end{equation}

For any complex place $v$ of $S$ (corresponding to  $\sigma: F \to \mathbb C$ and its conjugate), let $\{\gamma_1, \cdots, \gamma_{2g}\}$ be a basis of the free module $H_1(A_{\sigma}(\mathbb C), \mathbb Z)$ of rank $2g$ and define the period
\begin{equation}\label{grossp} p_{\sigma} (A, \eta) = \big | \textrm{det}( M_v(A, \eta)) \big |, \qquad  M_v(A, \eta) = \bigg(\bigg(\int_{\gamma_i}\omega_j, {\overline{{\int_{\gamma_i}{{\omega}_j}}}}~\bigg) \bigg)~;\end{equation} 
this period is non-zero and depends only on the differential form $\eta$ and the place $v$. 
The  period of $A$ (relative to $\eta$) is the real number
\begin{equation}\label{archperiod} 
P_{A, \infty}(\eta)  = \prod_{\sigma \in S_{\infty}} p_{\sigma}(A, \eta).
\end{equation}

\begin{remark}\label{wellknown} We recall the well known fact  that the determinant of $M_v(A, \eta)$ is $(\sqrt{-1})^{g}$ times a real number.   In the basis $\{\kappa_1, \cdots, \kappa_g\}$ of $H^0(A_{\sigma}(\mathbb C), \Omega^1)$ defined by $\int_{\gamma_j}\kappa_i = \delta_{ij}$ for $1\le i,j\le g$, the $2g\times 2g$-period matrix $(\int_{\gamma_i}\kappa_j, {\overline{{\int_{\gamma_i} \kappa_j}}})$ is a block matrix of the form
\[M'= \begin{bmatrix} I_g & I_g\\ \Omega & {\overline{\Omega}}
\end{bmatrix}, \qquad \Omega = \bigg(\int_{\gamma_{g+i}}\kappa_j\bigg)_{1\le i,j \le g} 
\]whose determinant is equal to the determinant of the $g \times g$-matrix ${\overline{\Omega}} - \Omega$ (all its entries have zero real part) and hence $\textrm{det}(M')$ is the product of a real number with $\sqrt{-1}^g$.  If $K$ is the change of basis matrix from the basis $\{w_1, \cdots, w_g\}$ to the basis $\{\kappa_1, \cdots, \kappa_g\}$, then
\[
M_v(A, \eta) = M'~\cdot~ \begin{bmatrix} K & 0\\ 0 & {\overline{K}}
\end{bmatrix} \implies \textrm{det} (M_v(A, \eta)) = \textrm{det}(M')~\cdot~\textrm{det}(K)~\cdot~{\overline{\textrm{det}(K)}}. \qed \]
\end{remark} 
\begin{lemma}\label{indep} The number 
\[P_{A, \infty}(\eta)~\cdot~ \mathbb N_{F/{\mathbb Q}}(\mathfrak a_{\eta}) =  \mathbb N_{F/{\mathbb Q}}(\mathfrak a_{\eta})~\cdot~\prod_{v \in S_{\infty}} p_{v}(A,\eta)\]
is independent of $\eta$. It will be denoted $P_{A, \infty}$. 
\end{lemma} 
\begin{proof} For any $0 \neq k \in F$ and any $v \in S_{\infty}$, one has
\[ p_{v}(A, k\eta)= \begin{cases} p_{v}(A, \eta)~\cdot~\sigma(k)& \mbox{if}~v\in S_{\mathbb R}\\ 
p_v(A, \eta)~\cdot~\overline{\sigma(k)}~\cdot~\sigma(k)&\mbox{if}~v\in S_{\mathbb C} \end{cases}.\]  
On the other hand, we see from (\ref{c-eta}) that $\mathfrak a_{k\eta} = (k)^{-1} \mathfrak a_{\eta}$. 
Since 
\begin{equation}\label{prod=norm} \prod_{\sigma \in \tilde{S}} \sigma(k) =  \mathbb N_{F/{\mathbb Q}}(k),\end{equation}
one has 
\[  \mathbb N_{F/{\mathbb Q}}(\mathfrak a_{k\eta})~\cdot~\prod_{\sigma \in S_{\infty}} p_{\sigma}(A,k\eta) =  \mathbb N_{F/{\mathbb Q}}(k)^{-1}~\cdot~ \mathbb N_{F/{\mathbb Q}}(\mathfrak a_{\eta}) ~\cdot~\prod_{\sigma \in S_{\infty}} p_{\sigma}(A,\eta) ~\cdot~(\prod_{\sigma \in \tilde{S}} \sigma(k))  = \mathbb N_{F/{\mathbb Q}}(\mathfrak a_{\eta}) ~\cdot~\prod_{\sigma \in {S}_{\infty} } p_{\sigma}(A,\eta).\]
\end{proof} 
\begin{definition} The global volume $P_A$ of $A$ is defined by
\[ P_A =  \frac{P_{A, fin}~\cdot~P_{A, \infty}}{|d_F|^{g/2}} = \frac{P_{A, fin}~\cdot~P_{A, \infty}(\eta)~\cdot~ \mathbb N_{F/{\mathbb Q}}(\mathfrak a_{\eta}) }{|d_F|^{g/2}}.\] 
 \end{definition} 
 \begin{remark}\label{conradbsd} It is known that $P_A$ is the volume $m_A(A(\mathbb A_F))$ of  the adelic points $A(\mathbb A_F)$ with respect to the (canonical) Tamagawa measure $m_A$.   If $A'$ is the Weil-restriction of $A$ to $F'\subset F$, then $L(A,s) = L(A',s)$. One can compare the terms in Conjecture \ref{bsd}  for $A$ and $A'$: 
\[ \Sha(A/F) \cong \Sha(A'/{F'}), \quad \Theta_{NT}(A) = \Theta_{NT}(A'), \quad A'(F') = A(F),\quad (A')^t (F') = A^t(F);\] 
see Remark \ref{allperiods} for the the comparison of $P_A$ and $P_{A'}$. Conjecture \ref{bsd} is compatible with restriction of scalars \cite{MR330174} \cite[\S 6 ]{Conrad}: Conjecture \ref{bsd} for $A$ over $F$ is equivalent to $A'$ over $F'$; see \cite[\S 6]{Conrad} for a beautiful exposition.\qed
 \end{remark} 
 
Let $A^t$ be the dual abelian variety. We write $\Theta_{NT}(A)$  for the determinant of the N\'eron-Tate height pairing
\begin{equation}\label{Nerontate} 
 \frac{A(F)}{A(F)_{tor}} \times \frac{A^t(F)}{A^t(F)_{tor}} \to \mathbb R.
 \end{equation} 
 Since this height pairing corresponds to the Poincar\'e divisor on $A \times A^t$ which is symmetric, it follows that $\Theta_{NT}(A) = \Theta_{NT}(A^t)$.  Following \S \ref{pairings}, one has \[ \Delta_{NT}(A) = \frac{\Theta_{NT}(A)}{([A(F)_{tor}]~\cdot~[A^t(F)_{tor}])}.\]
Finally, let $\Sha(A/F)$ denote the Tate-Shafarevich group of $A$ over $F$.

The BSD conjecture states that the order of vanishing of $L(A,s)$ at $s=1$ is equal to the rank of the Mordell-Weil group $A(F)$ of $A$. The strong BSD conjecture \cite{Tate1966} (we follow Gross's formulation \cite[Conjecture 2.10]{Gross82} which is shown by him to be equivalent to the one in \cite{Tate1966})  
concerns the special value $L^*(A,1)$ of  $L(A,s)$ at $s=1$. 
\begin{conjecture}\label{bsd}  $\Sha(A/F)$ is finite and  $L^*(A,1)$ satisfies 
\begin{equation}\label{bsd1} L^*(A,1) = \frac{P_A~\cdot~\Theta_{NT}(A)~\cdot~[\Sha(A/F)]}{([A(F)_{tor}]~\cdot~[A^t(F)_{tor}])} = P_A~\cdot~\Delta_{NT}(A)~\cdot~[\Sha(A/F)].\end{equation} 
 \end{conjecture}

\subsection{The case of everywhere good reduction}   If $A$ has good reduction everywhere, then  $c_v=1$ for all $v \in S$ and  $P_{A,fin} =1$; also, $A_v = \mathcal N_v = \mathcal N^0_v = B_v.$ Conjecture \ref{bsd} simplifies to 

\begin{conjecture}\label{bsdgood}  If $A$ has good reduction everywhere, then
\[L^*(A,1) = \frac{P_{A, \infty}(\eta) ~\cdot~ \mathbb N_{F/{\mathbb Q}}(\mathfrak a_{\eta})~\cdot~\Theta_{NT}(A)~\cdot~[\Sha(A/F)]}{{|d_F|^{g/2}}([A(F)_{tor}]~\cdot~[A^t(F)_{tor}])}.\]
\end{conjecture}

\section{Comparison of periods}\label{2periods}  In order to relate Conjectures \ref{main conjecture} and \ref{bsd}, it is necessary to compare the period $P_{J, \infty}$ of Lemma \ref{indep} with the determinant of $B(1,1)$ in (\ref{Bj1}).  Let $\omega_{\mathcal J}$ denote the projective $\mathcal O_F$-module of invariant differentials on  the N\'eron model $\mathcal J$ of the Jacobian $J$ of $X_0$.  Recall that the integral structure $i_{\omega_{\mathcal J}}$ (Definition \ref{im-integral})  uses the abelian group underlying $\omega_{\mathcal J}$ and forgets the $\mathcal O_F$-structure.
 \begin{definition} (i)  $det(\gamma_{\mathcal J}^*)$ is the determinant of the period isomorphism
 \begin{equation}\label{gamma-j*} \gamma^*_{\mathcal J}:~H^1_B(J_{\mathbb C}, \mathbb Z)^{+}_{\mathbb C} \xleftarrow{\sim} \omega_{\mathcal J} \otimes_{\mathbb Z} \mathbb C\end{equation} 
 calculated with respect to $H^1_B(J_{\mathbb C}, \mathbb Z)^{+}$ and the integral structure $i_{\omega_{\mathcal J}}$. 
 
\noindent (ii) $det_{\mathcal O_F}(\gamma_{\mathcal J}^*)$ is the determinant of (\ref{gamma-j*}) with respect to  $H^1_B(J_{\mathbb C}, \mathbb Z)^{+}$ and the $\mathcal O_F$-integral structure  $j_{ \omega_{\mathcal J} }$ (Definition \ref{jm-integral}).  
 \end{definition}
 The main result of this section is the following
 \begin{theorem}\label{allperiods} One has (up to a power of two)
 \[det(\gamma_{\mathcal J}^*)= \pm \frac{P_{J,\infty}}{|d_F|^{g/2}} =  \pm ({\sqrt{-1}})^{g~\cdot~r_2}~\cdot~\frac{P_{J, \infty}(\eta)~\cdot~ \mathbb N_{F/{\mathbb Q}}(\mathfrak a_{\eta}) }{\sqrt{d_F}^g}.\]
 \end{theorem} 
 \begin{remark}\label{allperiods-remark} (i) As we shall see below,  the factor $\sqrt{d_F}^g$ arises from the change of integral structures, the factor $\mathbb N_{F/{\mathbb Q}}(\mathfrak a_{\eta}) $ appears if $\omega_{\mathcal J}$ is not free, and  if $F$ is totally real, then $P_{J, \infty} = det_{\mathcal O_F}(\gamma^*_{\mathcal J})$ up to a power of two (because (\ref{grossreal}) - but not (\ref{gamma-j*}) - uses  $\pi_0(A_{\sigma}(\mathbb R))$, a finite two group). Finally, there is a factor of $\sqrt{-1}^{r_2~\cdot~g}$ because the definition of $P_{J,\infty}$ uses a different Betti lattice at complex places. 
 
(ii) As  $det(\gamma_{\mathcal J}^*)$ uses only the underlying abelian group of $\omega_{\mathcal J}$, one observes that the actual
$\mathcal O_F$-structure on $\omega_{\mathcal J}$ is irrelevant for Conjecture \ref{bsd}. This observation has the following implication. As in  Remark \ref{conradbsd}, let $A$ be an abelian variety over $F$ and let $A'$ its restriction to $F'\subset F$; let $\mathcal N$ and $\mathcal N'$ be their N\'eron models. Conjecture \ref{bsd} for an abelian variety $A$ over $F$ is equivalent to Conjecture \ref{bsd} for  $A'$ over $F'$.   A key step in the proof of the equivalence is showing that  \cite[\S 6]{Conrad} 
\[  \frac{P_{A,\infty}}{(|d_F|)^{g/2}} =  \frac{P_{A',\infty}}{(|d_{F'} |)^{\frac{g[F:F']}{2}}}.\]This is clear from Theorem \ref{allperiods} as $\omega_{\mathcal N}$ is isomorphic to $\omega_{\mathcal N'}$ as abelian groups and the period isomorphism $\gamma_{\mathcal N}$ is the same as $\gamma_{\mathcal N'}$. \qed \end{remark} 
 \begin{proof} (of Theorem \ref{allperiods}) 
 It is difficult to calculate or write down a formula directly for $det(\gamma_{\mathcal J}^*)$ as the integral structure $i_{\omega_{\mathcal J}}$ is not easy to describe explicitly.  So one has to use $j_{\omega_{\mathcal J}}$.  By Proposition \ref{euler-disc}, one has 
 \[ det_{\mathcal O_F}(\gamma^*_{\mathcal J}) = det(\gamma^*_{\mathcal J})~\cdot~\sqrt{d_F}^g.\] We shall now see how to compute $det_{\mathcal O_F}(\gamma_{\mathcal J})$ explicitly in the style of   \S \ref{arch-periodes}.

\subsubsection{The case that $N= \omega_{ \mathcal J}$ is free} If $N$ is free as a $\mathcal O_F$-module (its rank is $g$), then pick a $\mathcal O_F$-basis $\mathcal B = \{v_1, \cdots, v_g\}$ for it. If $\mathcal B'$ is a different basis for $N$, then the determinant $d(\mathcal B, \mathcal B')$ of the change of basis matrix is a unit in $\mathcal O_F$. 
Note that $\mathcal B$ is also a basis for the $F$-vector space $\omega_{\mathcal J}\otimes_{\mathcal O_F}F = H^0(J, \Omega^1)$. 
For any embedding $\sigma: F \to \mathbb C$, let  \[ \sigma(\mathcal B)= \{\sigma(v_1),\cdots, \sigma(v_g)\}\] be the image of $\mathcal B$ under the isomorphism (via $\sigma$) 
\[H^0(J, \Omega^1)\otimes_{F}\mathbb C \cong H^0(J_{\sigma}(\mathbb C), \Omega^1).\] 
So $\sigma(\mathcal B)$ is a basis for the complex vector space $H^0(J_{\sigma}(\mathbb C), \Omega^1)$. Just as a basis for $V$ and $W$ produces a basis for $V \times W$, we can combine the bases $\sigma (\mathcal B)$ for $\sigma \in \tilde{S}$ to get a basis $\tilde{\mathcal B}$ for the complex vector space
\[H^0(J_{\mathbb C}, \Omega^1) = H^0(J, \Omega^1)\otimes_{\mathbb Z} \mathbb C=  \prod_{\sigma \in \tilde{S}} H^0(J_{\sigma}(\mathbb C), \Omega^1).\]
The lattice ($\mathbb Z$-module) spanned by $\tilde{\mathcal B}$ represents the $\mathcal O_F$-integral structure $j_N$ on $V =N\otimes_{\mathbb Z}\mathbb C$ (see Definition \ref{jm-integral}). The determinant  $D(\mathcal B)$ of (\ref{gamma-j*}) computed with respect to the integral structure provided by $\tilde{\mathcal B}$ is equal to  $det_{\mathcal O_F}(\gamma^*_{\mathcal J})$. Let us show that $D(\mathcal B)$ is independent of the basis $\tilde{\mathcal B}$.

Given another basis $\mathcal B'$ for $N$, we have
\[D(\mathcal B')  = D(\mathcal B) \prod_{\sigma\in \tilde{S} }\sigma(d_{\mathcal B, \mathcal B'}).\] The determinant of the change of basis matrix from $\tilde{\mathcal B}$ to $\tilde{\mathcal B'}$ is 
\[\prod_{\sigma\in \tilde{S}}\sigma(d_{\mathcal B, \mathcal B'}) = \mathbb N_{F/{\mathbb Q}}d_{\mathcal B, \mathcal B'} =\pm 1\]
as $d_{\mathcal B, \mathcal B'}$ is a unit in $\mathcal O_F$.
Though the integral structure on  $N\otimes_{\mathbb Z} \mathbb C$ provided by $\tilde{\mathcal B}$ depends on the choice of $\mathcal B$, the number $D(\mathcal B)$ is independent of the basis $\mathcal{B}$:  Any two $\mathcal O_F$-bases of $N$ provide Euler-equivalent integral structures representing $j_N$ on $V=N\otimes_{\mathbb Z} \mathbb C$. Thus, in this case, for any basis $\mathcal B$, one has
\[ det(\gamma^*_{\mathcal J}) =  \frac{D(\mathcal B)}{{\sqrt{d_F}^g}}.\]

\subsubsection{The case that $N= \omega_{\mathcal J}$ is not free as a $\mathcal O_F$-module}  In this case, we proceed as in \S \ref{arch-periodes} using the top exterior power $\Lambda^g N$ of the projective $\mathcal O_F$-module $N$. Let $\mathcal B = \{v_1, \cdots, v_g\}$ be a $F$-basis of $H^0(J, \Omega^1)$ and put $\rho = v_1 \wedge \cdots \wedge v_g$. We have
\begin{equation}\label{brho1} 
\Lambda^g N = \mathfrak b_{\rho}~\cdot~ \rho \subset \Lambda^g H^0(J, \Omega^1)
\end{equation} 
where $\mathfrak b_{\rho}$ is a fractional ideal of $\mathcal O_F$. For any $k \neq 0 \in F$, one has an equality of fractional ideals
\begin{equation}\label{brho}
\mathfrak b_{\rho} = \mathfrak b_{k\rho}~\cdot~(k).
\end{equation} 

As before, let $\tilde{\mathcal B}$ denote the basis of $V=N\otimes_{\mathbb Z} \mathbb C$ obtained from $\mathcal B$. Let $D(\mathcal B)$ denote the determinant of (\ref{gamma-j*}) with respect to the basis $\tilde{\mathcal B}$ and $H^1_B(J_{\mathbb C}, \mathbb Z)^{+}$.  Using 
(\ref{brho1}), one has 
\[ det_{\mathcal O_F}(\gamma^*_{\mathcal J}) = D(\mathcal B)~\cdot~\mathbb N_{F/{\mathbb Q}} \mathfrak b_{\rho}.\]
Let us directly show that the right hand side is independent of the basis.  Given any basis $\mathcal B'=\{v_1',\cdots, v_g'\}$, then
 \[ D(\mathcal B') = D(\mathcal B) \prod_{\sigma \in \tilde{S} } \sigma(k) \qquad k:= d(\mathcal B, \mathcal B').\] 
 Note that $k= d(\mathcal B, \mathcal B')\in F^{\times}$ need not be a unit.

If $\rho'= v_1'\wedge\cdots\wedge v_g'$, then
\[\rho' = \rho~\cdot~k, \qquad k= d(\mathcal B, \mathcal B')\]
which gives
\[b_{\rho'} = b_{\rho}~\cdot~(k)^{-1}.\]
Thus, using (\ref{prod=norm}) as in Lemma \ref{indep}, we find 
\[D(\mathcal B')~\cdot~\mathbb N_{F/{\mathbb Q}} \mathfrak b_{\rho'} =  D(\mathcal B)~\cdot~( \prod_{\sigma \in \tilde{S} } \sigma(k))~\cdot~\mathbb N_{F/{\mathbb Q}} \mathfrak b_{\rho'} = D(\mathcal B)~\cdot~(\prod_{\sigma \in \tilde{S}} \sigma(k))~\cdot~\mathbb N_{F/{\mathbb Q}} \mathfrak b_{\rho}~\cdot~(\mathbb N_{F/{\mathbb Q}}(k))^{-1} =  D(\mathcal B)~\cdot~\mathbb N_{F/{\mathbb Q}} \mathfrak b_{\rho}.\] 
Thus, $D(\mathcal B)~\cdot~\mathbb N_{F/{\mathbb Q}} \mathfrak b_{\rho}$ is independent of $\mathcal B$. 

\subsubsection{Betti lattices} It should be clear that the above computation of $det_{\mathcal O_F}(\gamma^*_{\mathcal J})$ is exactly the computation of $P_{J,\infty}$ in \S \ref{arch-periodes}. In fact, for any $g$-dimensional abelian variety $A$ over $F$, one has (up to a power of two)
\begin{equation}\label{signs+} det_{\mathcal O_F}(\gamma^*_{\mathcal N}) = \pm ({\sqrt{-1}})^{g~\cdot~r_2}~\cdot~P_{A,\infty}.\end{equation}
 Here, as in \S \ref{arch-periodes}, $\mathcal N$ is the N\'eron model of $A$ and $\gamma^*_{\mathcal N}$ is the period isomorphism 
\[\gamma^*_{\mathcal N}:~H^1_B(A_{\mathbb C}, \mathbb Z)^{+}_{\mathbb C} \xleftarrow{\sim} \omega_{\mathcal N} \otimes_{\mathbb Z} \mathbb C.\]
 Let $B= \{w_1, \cdots, w_g\}$ be a $F$-basis of $H^0(A, \Omega^1)$  and put $\eta = w_1\wedge\cdots \wedge w_g$. The discrepancy in (\ref{signs+}) arises from the fact that the Betti lattices are different.

For a real place of $S$ corresponding to $\sigma: F \to \mathbb R$,  consider the period isomorphism
\[ \gamma_{\sigma}= H^1_B(A_{\sigma}(\mathbb C) , \mathbb Z)^{+}_{\mathbb C} \xleftarrow{\sim} H^0(A_{\sigma}, \Omega^1).\]
For a complex place $v$ of $S$ (corresponding to $\sigma: F \to \mathbb C$ and its conjugate $c\sigma$), consider the complex abelian variety 
 \[A_v = A_{\sigma} \times A_{c\sigma}\]and the period isomorphism
 \[\gamma_{v}:  H^1_B(A_v, \mathbb Z)^+_{\mathbb C} \xleftarrow{\sim} H^0(A_v, \Omega^1).\]
 
 For real places, the groups in $\gamma_{\sigma}$ are the same as the one used in  (\ref{grossreal});  in this case, $p_{\sigma}(A, \eta)$ is the determinant\footnote{up to a power of two} of $\gamma_{\sigma}$ relative to the basis $H^1_B(A_{\sigma}(\mathbb C) , \mathbb Z)^{+}$ and $B=  \{w_1, \cdots, w_g\}$.  
 
 For complex places $v$,  (\ref{grossp}) uses $H^1(A_{\sigma}, \mathbb Z)$ instead of $H^1_B(A_v, \mathbb Z)^+$. In this case, $(\sqrt{-1})^g~\cdot~p_{v}(A, \eta)$ is the determinant of $\gamma_v$ relative to the basis $H^1_B(A_v, \mathbb Z)^+$ and $B$. This is what leads to the discrepancy between the two invariants. 
  What follows is presumably well known, but we include it for sake of completeness. 
  
   Complex conjugation \[c: A_{\sigma} \to A_{c\sigma}\] simply permutes the factors of ${A}_v$. One has the decomposition
\[H_1({A_v}(\mathbb C), \mathbb Z) = H_1(A_{\sigma}(\mathbb C), \mathbb Z) \oplus H_1(A_{c\sigma} (\mathbb C), \mathbb Z), \qquad H^0({A_v}, \Omega^1) = H^0(A_{\sigma}, \Omega^1) \oplus H^0(A_{c\sigma}, \Omega^1).\] Concretely, a basis for  $H_1({A_v}, \mathbb Z)$ is given by
\[\{(\gamma_1, 0), (\gamma_2, 0), \cdots, (\gamma_{2g}, 0), (0, c\gamma_1), (0,c\gamma_2), \cdots, (0,c\gamma_{2g})\};\]
a basis for $H^0(A_v, \Omega^1)$ is given by (here the image of $w_j$ is denoted as $w_j$ on $A_{\sigma}$ and as $cw_j$ on $A_{c\sigma}$) 
\[\{(w_1, 0), (w_2, 0), \cdots, (w_g, 0), (0, cw_1), (0,cw_2), \cdots, (0, cw_g)\}.\]
So a basis for $H_1(A_v, \mathbb Z)^{+}$  is given by 
\[\{(\gamma_1, c\gamma_1), \cdots, (\gamma_{2g},  c\gamma_{2g})\}.\]
For the complex place $v$ and $\eta$, the determinant of $\gamma_v$ relative to the basis $H^1_B(A_v, \mathbb Z)^+$ and $B$ is the determinant of the $2g \times 2g$-matrix $M$
\[
M=\begin{bmatrix} \int_{(\gamma_1, c\gamma_1)} (w_1, 0) &  \int_{(\gamma_1, c\gamma_1)} (w_2, 0) & \cdots & \int_{(\gamma_1, c\gamma_1)} (w_g, 0) & \int_{(\gamma_1, c\gamma_1)} (0, cw_1) & \cdots & \int_{(\gamma_1, c\gamma_1)} (0, cw_g)\\
\int_{(\gamma_2, c\gamma_2)} (w_1, 0) &  \int_{(\gamma_2, c\gamma_2)} (w_2, 0) & \cdots & \int_{(\gamma_2, c\gamma_2)} (w_g, 0) & \int_{(\gamma_2, c\gamma_2)} (0, cw_1) & \cdots & \int_{(\gamma_2, c\gamma_2)} (0, cw_g)\\
\cdots &\cdots  & \cdots &\cdots &\cdots &\cdots & \cdots\\
\cdots &\cdots  & \cdots &\cdots &\cdots &\cdots & \cdots\\
\cdots &\cdots  & \cdots &\cdots &\cdots &\cdots & \cdots\\
\int_{(\gamma_{2g}, c\gamma_{2g})} (w_1, 0) &  \int_{(\gamma_{2g}, c\gamma_{2g})} (w_2, 0) & \cdots & \int_{(\gamma_{2g}, c\gamma_{2g})} (w_g, 0) & \int_{(\gamma_{2g}, c\gamma_{2g})} (0, cw_1) & \cdots & \int_{(\gamma_{2g}, c\gamma_{2g})} (0, cw_g)
\end{bmatrix} 
\]
Observe that (in a product $V \times W$ of manifolds, a differential form from $W$ does not pair with a cycle from $V$)
\[\int_{(\gamma_i, c\gamma_i)} (w_j, 0) = \int_{\gamma_i} w_j, \qquad \int_{(\gamma_i, c\gamma_i)} (0, cw_j) = \int_{c\gamma_i} c w_j.\]

So the matrix $M$ can be rewritten as
\[
M=\begin{bmatrix} \int_{\gamma_1} w_1 &  \int_{\gamma_1} w_2 & \cdots & \int_{\gamma_1} w_g &  \int_{c\gamma_1} cw_1 & \cdots & \int_{ c\gamma_1} cw_g\\
\int_{\gamma_2} w_1 &  \int_{\gamma_2} w_2 & \cdots & \int_{\gamma_2} w_g &\int_{c\gamma_2} cw_1 & \cdots & \int_{c\gamma_2} cw_g\\
\cdots &\cdots  & \cdots &\cdots &\cdots &\cdots & \cdots\\
\cdots &\cdots  & \cdots &\cdots &\cdots &\cdots & \cdots\\
\cdots &\cdots  & \cdots &\cdots &\cdots &\cdots & \cdots\\
\int_{\gamma_{2g}} w_1 &  \int_{\gamma_{2g}} w_2 & \cdots & \int_{\gamma_{2g}} w_g & \int_{ c\gamma_{2g} }  cw_1 & \cdots & \int_{ c\gamma_{2g}} cw_g
\end{bmatrix}. 
\]
As $A_{\sigma}$ and $A_{c\sigma}$ are complex conjugate abelian varieties, one has
\[\int_{ c\gamma_{i} }  cw_j = \overline{\int_{\gamma_{i} }  w_j};\]
it follows that $M$ and $M_v(A, \eta)$ of (\ref{grossp}) satisfy 
\[ \textrm{det} (M) =  \textrm{det}(M_v(A, \eta)).\]
Remark \ref{wellknown} shows that  
 \[\prod_{v\in S_{\infty}} \textrm{det}(M_v (A, \eta))\]
 is a real number multiplied by $({\sqrt{-1}})^{g~\cdot~r_2}$.   This proves the identity (\ref{signs+}); combining it with Proposition \ref{euler-disc}  provides (up to a power of two)
 \[ det(\gamma^*_{\mathcal J}) = \frac{ det_{\mathcal O_F}(\gamma^*_{\mathcal J})}{\sqrt{d_F}^g} = \pm ({\sqrt{-1}})^{g~\cdot~r_2} \frac{P_{J,\infty}}{\sqrt{d_F}^g}.\]
 Theorem \ref{allperiods} follows using the well known result that the sign of $d_F$ is $(-1)^{r_2}$. \end{proof}

\begin{proposition}\label{yetanother} (i) $det_{\mathcal O_F}(\gamma_{\mathcal J}^*) $ is equal to the determinant $det_{\mathcal O_F}(\gamma_{\mathcal J})$ of 
\[\gamma_{\mathcal J}:~H^1_B(J_{\mathbb C}, \mathbb Z(1))^{+}_{\mathbb C} \xrightarrow{\sim} \textrm{Lie}~(\mathcal J) \otimes_{\mathbb Z} \mathbb C\]
calculated using the $\mathcal O_F$-integral structure $j_{\textrm{Lie}~(\mathcal J)}$ (Definition \ref{jm-integral}) and the lattice
$H^1_B(J_{\mathbb C}, \mathbb Z(1))^{+}$.  

(ii) $det(\gamma_{\mathcal J}^*)= det(\gamma_{\mathcal J})$. 
\end{proposition} 
\begin{proof} (i) For any complex abelian variety $A$, the dual abelian variety $A^t$ satisfies
\begin{align*} H_1(A^t, \mathbb Z) = H^1_B(A, \mathbb Z(1))& & \textrm{Lie}~A = H^1(A^t, \mathcal O) \\
\textrm{Hom}(H^1_B(A^t, \mathbb Z), \mathbb Z)  = H^1_B(A, \mathbb Z(1)) && \textrm{Hom}_{\mathbb C}(H^1(A, \mathcal O), \mathbb C) = H^0(A^t, \Omega^1).
\end{align*} 
This shows that the period isomorphism 
\[\gamma_A: H^1_B(A, \mathbb Z(1))^+_{\mathbb C} \xrightarrow{\sim} H^1(A, \mathcal O)\]
is dual to the map
\[\gamma^*_{A^t}: H^1_B(A^t, \mathbb Z)^+_{\mathbb C} \xleftarrow{\sim} H^0(A^t, \Omega^1).\]
Applying these to the self-dual Jacobian $J_{\mathbb C}$ gives that the lattice $H^1_B(J_{\mathbb C}, \mathbb Z(1))^+$ is dual to the 
lattice $H^1_B(J_{\mathbb C}, \mathbb Z)^+$. The natural duality of the projective $\mathcal O_F$-modules 
$\textrm{Lie}~(\mathcal J)$ and $\omega_{\mathcal J}$ shows, using Proposition \ref{euler-disc}, that the 
 $\mathcal O_F$-integral structures $j_{\textrm{Lie}~(\mathcal J)}$ and $j_{ \omega_{\mathcal J} }$ are dual. As the determinants of dual 
 maps computed with respect to dual lattices are equal, the result follows. 
 
 (ii) Proposition \ref{euler-disc} and (i) imply $det(\gamma_{\mathcal J}^*)~\cdot~(\sqrt{d_F})^g = det_{\mathcal O_F}(\gamma_{\mathcal J}^*) =det_{\mathcal O_F}(\gamma_{\mathcal J})  =det(\gamma_{\mathcal J}) ~\cdot~(\sqrt{d_F})^g$.
\end{proof}

\section{Proof of Theorem \ref{Main theorem} in a special case}\label{Proof} 
Throughout this section, we assume that $\pi:X \to S$ is smooth and $X_0(F)$ is non-empty. We prove Theorem \ref{Main theorem} in this case by computing $H^*_W(X, \mathbb Z(1))$ and  then use it to compare $\chi(X,1)$ with Conjecture \ref{bsdgood} for the Jacobian $J$ of $X_0$.

\subsection{The groups $H^*_W(X, \mathbb Z(1))$} By \S \ref{Euler1}, we need to understand $\textrm{Pic}(X)$ and $\textrm{Br}(X)$.  

Any $x\in X_0(F)$ provides $ {\rm Pic}(X_0) = J(F) \times \mathbb Z$.

\begin{proposition}\label{sha=Br} One has (neglecting two-torsion) 
\[ (i) ~\textrm{Pic}(X) \cong \mathbb Z \times J(F) \times \textrm{Pic}(S),\quad [\textrm{Pic}(X)_{tor}] = h~\cdot~[J(F)_{tor}], \qquad (ii)~{\rm Br}(X) \xrightarrow{\sim} \Sha(J/F).\]

\end{proposition} 
\begin{proof} As $\pi:X\to S$ is smooth proper, any $x\in X_0(F)$ provides a splitting of the map $\textrm{Pic}(S) \to \textrm{Pic}(X)$.  The identity component $\textrm{Pic}^0_{X/S}$ of the relative Picard scheme $\textrm{Pic}_{X/S}$ is the N\'eron model $\mathcal J$ of $J$ by  \cite[Theorem 1, page 264]{BLR} and 
$\textrm{Pic}_{X/S}(S) = \textrm{Pic}^0_{X/S}(S) \times \mathbb Z$ by \cite[Theorem 1, page 252]{BLR}.  As $\pi_*\mathcal O_X =\mathcal O_S$, \cite[Proposition 4, page 204]{BLR} says $\textrm{Pic}_{X/S}(S) = {\rm Pic}(X)/{{\rm Pic}(S)}$ and 
$\textrm{Pic}_{X/S}(F) = {\rm Pic}(X_0)$. This proves (i).   
For (ii), we note that \cite[Theorem 3.1]{Tate1966} provides an exact (modulo two-torsion) sequence
\[ 0 \to {\rm Br}(S) \to {\rm Br}(X) \to \Sha(J/F) \to 0;\] by Class Field Theory, $\textrm{Br}(S)$ is a finite $2$-group. \end{proof} 

\subsection{The Euler characteristic $\chi_{B}(X,1)$}\label{B11}We now compute the Euler characteristic of (\ref{Bj1}) 
\[\chi_{B}(X,1) = \frac{\chi(B(0,1)) \chi(B(2,1))}{\chi(B(1,1))}.\]

\begin{remark} As $\pi:X\to S$ is smooth, $H^1(X_s, \mathcal O_{X_s})$ has dimension $g$ over the residue field $k(s)$ for all points $s\in S$. Since $S$ is a reduced noetherian scheme, $\pi: X\to S$ is proper, and $\mathcal O_X$ is a coherent sheaf on $X$ flat over $S$, a result of Grothendieck (see L.~Illusie's article \cite[8.5.18]{FGA}) states that $R^1\pi_*\mathcal O_X $ is locally free and that $R^1\pi_*\mathcal O_X\otimes k(s) \cong H^1(X_s, \mathcal O_{X_s})$.  So the $\mathcal O_F$-module $H^1(X, \mathcal O)$ is projective of rank $g$.\qed
\end{remark} 

\subsubsection{The Euler characteristic of $B(1,1)$}  

We write $det(\gamma_X)$ for the determinant of 
\begin{equation}\label{gammaX} \gamma_X: H^1_B(X_{\mathbb C}, \mathbb Z(1))^{+}_{\mathbb C} \xrightarrow{\sim} H^1(X, \mathcal O)\otimes_{\mathbb Z} \mathbb C.
\end{equation}  with respect to the lattices $H^1_B(X_{\mathbb C}, \mathbb Z(1))^{+}$ and the abelian group underlying $H^1(X, \mathcal O)$. 

\begin{lemma}\label{jduality} The projective $\mathcal O_F$-modules $H^1(\mathcal J, \mathcal O)$ and $\omega_{\mathcal J}$ are dual. 
\end{lemma} 
\begin{proof} As $\mathcal J \to S$ is an abelian scheme, its relative Picard scheme $\textrm{Pic}_{\mathcal J/S}$ exists  \cite[Theorem 5, p. 234]{BLR} and its identity component $\textrm{Pic}^0_{\mathcal J/S}$ is the dual abelian scheme $\mathcal J^t \to S$. The $S$-points $\textrm{Lie}~\textrm{Pic}_{\mathcal J}$ of its Lie algebra\footnote{For any group scheme $G$ over $S$, we write $\textrm{Lie}~G$ for the $S$-points $Lie~G(S)$ of the Lie algebra $Lie~G$ over $S$.}  ${Lie}~\textrm{Pic}_{\mathcal J}$  satisfies \cite[Theorem 1, p.~231]{BLR}, \cite[Proposition 1.1 (d)]{MR2092767}
\begin{equation} \textrm{Lie}~(\mathcal J^t)= \textrm{Lie~(Pic}_{\mathcal J/S}) \xrightarrow{\sim}  H^1(\mathcal J, \mathcal O).
\end{equation}

Since $J^t$ is the generic fiber of $\mathcal J^t$, the self-duality $J \cong J^t$ shows $\mathcal J \cong \mathcal J^t$ and $\textrm{Lie}~(\mathcal J^t) \cong  \textrm{Lie}~(\mathcal J)$.  Combining this with the  natural duality between $ \textrm{Lie}~(\mathcal J)$ and $\omega_{\mathcal J}$ \cite[Proposition 1.1 (c)]{MR2092767} proves the lemma.
\end{proof}

\begin{proposition} \label{bone1} One has
\[\chi(B(1,1)) =  \frac{P_{J, fin}}{P_{J}}= \frac{{\sqrt{|d_F|}}^g}{P_{J,\infty}(\eta)~ \mathbb N_{F/{\mathbb Q}}(\mathfrak a_{\eta})}.\]
\end{proposition} 
\begin{proof}  Any $x\in X_0(F)$  gives a map
\[ \beta_x: X \to \mathcal J, \qquad X_0 \to J;\]the induced map $\beta$ on cohomology is independent of the choice of $x$; for example, it provides the following isomorphism
\[\beta: H^1(\mathcal J, \mathcal O) \xrightarrow{\sim} H^1(X, \mathcal O)\]
of $\mathcal O_F$-modules.    
which fits into a commutative diagram (using Lemma \ref{jduality}) 
  \begin{equation}
    \label{eq:1}
    \begin{tikzcd}
   H^1_B(J_{\mathbb C}, \mathbb Z(1))^{+}_{\mathbb C} \ar[r,"\gamma_{\mathcal J}"] \ar[d,"\beta"'] &
     H^1(\mathcal J, \mathcal O)\otimes_{\mathbb Z} \mathbb C  \ar[d,"\beta"] \\
      H^1_B(X_{\mathbb C}, \mathbb Z(1))^{+}_{\mathbb C} \ar[r,"\gamma_X"']        & H^1(X, \mathcal O)\otimes_{\mathbb Z} \mathbb C;
    \end{tikzcd}
  \end{equation}
as the vertical maps $\beta$ are isomorphisms, the integral structure corresponding to the map $\gamma_{\mathcal J}$ (top row) is isomorphic to the one corresponding to $\gamma_X$ (bottom row). As they are torsion-free, Theorem \ref{allperiods} and Proposition \ref{yetanother} show  that \[\chi(\gamma_X) = \chi(\gamma_{\mathcal J}) =  det(\gamma_{\mathcal J}) = \frac{P_{J, \infty}}{|d_F|^{g/2}}.\]

Our convention for $B(1,1)$ is that $\textrm{Ker}(\gamma_M)$ is in even degree, say degree zero. We have 
\[B(1,1):~\underset{\textrm{degree one}}{H^1_B(X_{\mathbb C}, \mathbb Z(1))^{+}_{\mathbb C}} \to \underset{\textrm{degree two}}{H^1(X, \mathcal O)\otimes_{\mathbb Z} \mathbb C}, \qquad det (B(1,1))  = \frac{{\sqrt{|d_F|}}^g}{P_{J,\infty}(\eta)~ \mathbb N_{F/{\mathbb Q}}(\mathfrak a_{\eta})} = \chi(B(1,1)). \]
\end{proof}

Combining Proposition \ref{bone1} with (\ref{bzero1}) and (\ref{btwo1}) yields the 
\begin{proposition}\label{B-Euler} We have
\[\frac{1}{\chi_B(X,1)}  =  \frac{\chi(B(1,1))}{\chi(B(0,1)) \chi(B(2,1))} =\frac{(2\pi i)^{r_2}}{\sqrt{d_F}}~\cdot~ \frac{{\sqrt{|d_F|}}^g}{P_{J,\infty}(\eta)~ \mathbb N_{F/{\mathbb Q}}(\mathfrak a_{\eta})} .\]
\end{proposition}

\subsection{The determinant of the sequence $C(1)$}\label{c11} 
Using (\ref{weil2}), (\ref{weil3}), (\ref{weil4}),  and Proposition \ref{sha=Br}, the torsion in (\ref{C1}) satisfies
\[\chi(C(1)_{tor}) = \frac{w~\cdot~[Br(X)]}{[\textrm{Pic}(X)_{tor}]^2 } = \frac{w~\cdot~[Br(X)]} {h~\cdot~[J(F)_{tor}]~\cdot~h~\cdot~[J(F)_{tor}]}.\]
As  $\textrm{Hom}(\textrm{Pic}(X), \mathbb Z)\otimes\mathbb C \cong \mathbb C \times \textrm{Hom}(J(F), \mathbb Z)\otimes\mathbb C$ by Proposition \ref{sha=Br}, 
 the sequence $C(1)_{\mathbb C}$ of (\ref{C1}) 
\begin{equation}\label{c1new} 
0 \to  \underset{degree~zero}{\mathcal O_F^{\times}\otimes\mathbb C} \to \underset{degree~one}{\textrm{Coker}(\gamma_{M_0})} \to \underset{degree~two}{\textrm{Hom}(\textrm{Pic}(X), \mathbb C)} \xrightarrow{h} \underset{degree~three}{{\textrm{Pic}(X)}\otimes\mathbb C} \to \underset{degree~four}{\textrm{Ker}(\gamma_{M_2})} \to \underset{degree~five}{\textrm{Hom}(\mathcal O_F^{\times}, \mathbb C)} \to 0.
\end{equation} 
breaks up into exact sequences
\begin{eqnarray*} 0 \to  \underset{degree~zero}{\mathcal O_F^{\times}\otimes\mathbb C} \to \underset{degree~one}{\textrm{Coker}(\gamma_{M_0})} \to \mathbb C \to 0,\\
\underset{degree~two}{\textrm{Hom}(\textrm{Pic}^0(X), \mathbb C)} \xrightarrow{h} \underset{degree~three}{{\textrm{Pic}^0(X)}\otimes\mathbb C},\\
0 \to \mathbb C \to \underset{degree~four}{\textrm{Ker}(\gamma_{M_2})} \to \underset{degree~five}{\textrm{Hom}(\mathcal O_F^{\times}, \mathbb Z)\otimes\mathbb C} \to 0.
\end{eqnarray*}
We compute $det(C(1))$  in three steps:  
\begin{itemize} 
 \item by  (\ref{gamma-m0}), the first sequence is  (\ref{conj-S-1}) and so has determinant $R$.

\item By \cite[Conjecture 2.3.7]{SL2020} or (\ref{Cr}), the map $h$ is the Arakelov intersection pairing \[\textrm{Hom}(J(F), \mathbb Z)\otimes{\mathbb C}  \to J(F)\otimes{\mathbb C}. \]
By Faltings-Hriljac \cite{MR740897, MR778087}, $-h$ is the N\'eron-Tate height pairing (\ref{Nerontate}) and so $det(h) = \pm \Theta_{NT}(J)$.  

\item As $t_{M_2} = H^2(X_0, \mathcal O) =0$, one has $\textrm{Ker}(\gamma_{M_2}) = H^2(X_{\mathbb C}, \mathbb C(1))^{+} \cong \mathbb C^{r_1+r_2}$. 
Using (\ref{weil4}), the third sequence is (\ref{conj-S-zero})  and so has determinant $R$.

\end{itemize}

\begin{proposition}\label{C-Euler} We have
\[ det(C(1)) = \frac{R \times R}{\Theta_{NT}} \qquad \chi (C(1)) = \frac{det(C(1))}{\chi(C(1)_{tor})}= \frac{h \times h  [J(F)_{tor}] [J(F)_{tor}]}{[Br(X)] w } \frac{R \times R}{\Theta_{NT}} .\]
\end{proposition}

\subsection{Completion of the proof of the main theorem in a special case}\label{zbsd} In this case, $J$ has good reduction everywhere as $\textrm{Pic}^0_{X/S}$ is the N\'eron model of $J$; the group scheme $\Phi_v$ is trivial, the Tamagawa numbers $c_v =1$ for all $v$ and $P_{J,fin} =1$.      
 Conjecture \ref{bsd} for $L(J,s)$ becomes Conjecture \ref{bsdgood}: 
\[L^*(J,1) = \frac{P_{J,\infty}(\eta)~ \mathbb N_{F/{\mathbb Q}}(\mathfrak a_{\eta}) }{(\sqrt{|d_F|})^{g}} \frac{[\Sha(J/F)]~\cdot~\Theta_{NT}}{[J(F)_{tor}] [J(F)_{tor}]}.\]  

We recall that
\[ \zeta(X,s) = \prod_{v\in S}\zeta(X_v, s), \qquad \zeta(X_v, s) = \frac{P_1(X_v,s)}{P_0(X_v, s) P_2(X_v, s)}. \]
As $\pi:X \to S$ is smooth and proper, one has 
\[P_0(X_v,s) = (1-q_v^{-s}), \quad P_1(C_v, t) = P_1(J_v, t),\quad P_2(X_v,s) = (1-q_v^{1-s})\]and 
\begin{equation}\label{onemore} 
\prod_{v\in S}\frac{1}{P_0(X_v,s)} = \zeta(S,s), \qquad \prod_{v\in S}\frac{1}{P_2(X_v,s)} = \zeta(S,s-1),  \quad L(J, s) = \prod_{v\in S} \frac{1}{P_1(X_v,s)}.\end{equation} 
Using (\ref{cnf0}), (\ref{cnf1}), we see that Conjecture \ref{bsdgood} for $J/F$ is equivalent to the equality
\[ \zeta^*(X,1) = \frac{\zeta^*(S,1) \zeta^*(S,0)}{L^*(J, 1)} =   {\frac{2^{r_1} (2\pi)^{r_2}hR}{w \sqrt{|d_F|}}}~\cdot~  {\frac{hR}{w}}~\cdot~{\frac{ (\sqrt{|d_F|})^{g} [J(F)_{tor}] [J(F)_{tor}]}{[\Sha(J/F)] \Theta_{NT} P_{J,\infty}(\eta)~ \mathbb N_{F/{\mathbb Q}}(\mathfrak a_{\eta}) }}.\]
 
 Conjecture \ref{main conjecture}  for $\zeta^*(X,1)$ says
 \[ \zeta^*(X,1) = \chi(X,1) = \frac{\chi_{A,C}(X,1)}{ \chi_B(X,1)} = \frac{{\chi(C(1))}}{\chi_B(X,1)~\cdot~ \chi(A'(3,1))}.\]

  Propositions \ref{B-Euler} and \ref{C-Euler} show
\[\chi(X,1) =  {\frac{1}{w}}~\cdot~{\frac{[\textrm{Pic}(X)_{tor}]^2}{w ~[\textrm{Br}(X)]}~\frac{R \times R}{\Theta_{NT}}}~\cdot~  {\frac{(2\pi i )^{r_2}}{\sqrt{d_F}}  \frac{{\sqrt{|d_F|}}^g}{P_{J,\infty}(\eta)~ \mathbb N_{F/{\mathbb Q}}(\mathfrak a_{\eta}) }}. \]

Proposition \ref{sha=Br} shows that Conjecture \ref{main conjecture} is equivalent to Conjecture \ref{bsdgood} for $J$ over $F$, up to powers of two.  This proves Theorem \ref{Main theorem} in the special case $\pi:X\to S$ is smooth and $X_0(F)$ is non-empty. \qed

\section{Towards the main theorem in the general case: arbitrary reduction}\label{tools+}  The assumption $\pi:X \to S$ is smooth and $X_0(F)$ is non-empty leads to several simplifications such as $H^1(X, \mathcal O)$ is torsion-free,  $\textrm{Br}(X) \cong \Sha(J/F)$, the intersection pairing on $\textrm{Pic}^0(X)$ is the N\'eron-Tate pairing on $J(F)$, and $\textrm{Pic}^0(X_0) \cong J(F)$. The failure of these identities complicate the proof of Theorem \ref{Main theorem} in the general case. In order to handle these complications, we now recall certain results from \cite{MR528839, MR2092767, MR2125738, MR3858404, Geisser1, MR3586808, MR1917232}.  

\subsection{Index and Period} Let $C$ be a smooth proper curve over a field $K$ and $T=$~Spec~$K$.

\begin{definition} 

(i) The index $\delta_K$ of $C$ over $K$ is the least positive degree of a $K$-rational divisor on $C$. If 
\[d: \textrm{Pic}(C) \to \mathbb Z\]is the degree map, then $\delta_K = [\textrm{Coker}(d)]$. The kernel of $d$ is $\textrm{Pic}^0(C)$.

(ii) The period $\delta'_K$ of $C$ over $K$ is the least positive degree of a $K$-rational divisor class on $C$.  \end{definition} 

\subsection{Curves over local and global fields}  \cite{MR1717533, MR2092767}

Let $K$ be a $p$-adic field with ring of integers $\mathcal O_K$ and residue field $k(v) = \mathbb F_q$. Write $T = {\rm Spec}~\mathcal O_K$ and $v: {\rm Spec}~\mathbb F_q \to T$ for the closed point of $T$.  Let $f: C \to T$ be a flat projective morphism with $C$ regular, $f_*\mathcal O_C = \mathcal O_T$ and $C_K \to {\rm Spec}~K$ a geometrically connected smooth curve of genus $g>0$.  Let $\mathcal N$ be the N\'eron model of the Jacobian $J$ of $C_K$. We write $\Phi_v=\pi_0(\mathcal N_v)$ for the finite group scheme of connected components of the special fiber $\mathcal N_v$ over $v$ and $c_v $ is the order of $\Phi_v(k(v))$. 

Consider the  map $Lie~(\phi): R^1f_* \mathcal O_C \to Lie~\mathcal N$ of coherent sheaves on $T$ and  the induced map on $T$-points 
\[{\rm Lie}(\phi): H^1(C, \mathcal O_C) \to {\rm Lie}~(\mathcal N).\]We will need the following result in \S \ref{Proof2}:
\begin{theorem}\label{raynaud1} \cite[Theorem 3.1]{MR2092767} 
(i) the kernel of  ${\rm Lie}~(\phi)$ is the torsion subgroup of $H^1(C, \mathcal O_C)$. 

(ii) the kernel and cokernel of $Lie~(\phi)$ are torsion sheaves on $T$ of the same length.  
\end{theorem}

Let $\Gamma_i$ ($i\in I$) be the irreducible components of the special fiber $C_v$. 
\begin{itemize}
\item $d_i$ is the multiplicity of $\Gamma_i$ in $C_v$.
\item $e_i$ is the geometric multiplicity of $\Gamma_i$ in $C_v$.
\item $r_i$ is the number of irreducible components of $\Gamma_i \times \textrm{Spec}~\overline{\mathbb F_q}$.
\end{itemize} 
There are canonical maps  \cite[Proposition 1.9]{MR1717533} $\alpha_C: \mathbb Z^I \to \mathbb Z^I$ (defined using intersection on $C$) and $\beta_C: \mathbb Z^I \to \mathbb Z$. Let $d$ be the gcd of the set $\{d_i, i\in I\}$ and $d'$ be the gcd of the set $\{r_id_i, i\in I\}$. 
\begin{theorem}\label{raynaudll} \cite[Theorem 1.7, Corollary 1.12]{MR1717533} (i) There is an exact sequence
\begin{equation*} 0 \to \frac{\textrm{Ker}~( \beta_C)}{\textrm{Im}~(\alpha_C)} \xrightarrow{h} \Phi_v(\mathbb F_q) \to \frac{cd\mathbb Z}{d'\mathbb Z} \to 0;\end{equation*}
here $c=1$ if $d'$ divides $g-1$ and $c=2$ otherwise. 

(ii) if $C_K(K)$ is not empty, then $d'=d$ and $h$ is an isomorphism. 
 \end{theorem}

It is known that  \cite[\S 9.1, Theorem 1.23]{MR1917232} that the intersection pairing on $R_v$ 
\[R_v:= \frac{\mathbb Z^I}{\mathbb Z} = \textrm{Coker}~(\mathbb Z \to \mathbb Z^I) \qquad 1\mapsto \sum_{i\in I} d_i \Gamma_i\]
is negative definite. The following corollary of Theorem \ref{raynaudll} and \cite[Theorem 1.11]{MR1717533} on $\Delta(R_v)$ (defined as in \S \ref{pairings}) is due to  Flach-Siebel. 
\begin{corollary}\label{fls} (Flach-Siebel \cite[Lemma 17]{flachsiebel}) If  $\delta_v$ is the index of $C_K$ and $\delta'_v$ is the period of $C_K$, then 
\[\Delta(R_v) = \frac{c_v}{\delta_v~\cdot~\delta'_v}  \prod_{i\in I} r_i.\]
In terms of the Arakelov intersection pairing on $R_v$ (see \cite[p.~390]{MR740897} or \cite[(3.7)]{MR861983}), one has
\[\Delta_{ar}(R_v) = \Delta(R_v) (\log q_v)^{\#G_v-1}.\] 
\end{corollary} 

\begin{lemma}\label{cohflat} Let $L_v(J, t)$ be the local $L$-factor of $J$ as in (\ref{localfactor}).
The zeta function 
\[Z(C_v,t) =  \frac{P_1(C_v,t)}{P_0(C_v,t)~\cdot~P_2(C_v,t)}\]  of $C_v$ satisfies
\[Z(C_v,t) = \frac{L_v(J,t)}{(1-t)~\cdot~\prod_{i\in I}(1-{(qt)}^{r_i})}.\] 
\end{lemma} 

\begin{proof} As $C_K$ is smooth projective and geometrically connected, $C_v$ is geometrically connected. So $P_0(C_v,t) = (1-t)$. 

We next show that $P_1(C_v,t)$  equals  $L_v(J,t)$ of (\ref{localfactor}). For any prime $\ell$ coprime to $q$, the Kummer sequence and the perfect pairing (Poincar\'e duality) 
\[H^1_{et}(C_K \times \bar{K}, \mathbb Q_{\ell}) \times H^1_{et}(C_K \times \bar{K}, \mathbb Q_{\ell}(1) ) \xrightarrow{\cup} H^2_{et}(C_K \times \bar{K}, \mathbb Q_{\ell}(1) ) \cong \mathbb Q_{\ell}\] 
provide isomorphisms of $\textrm{Gal}(\bar{K}/K)$-representations: 
\[H^1_{et}(C_K \times \bar{K}, \mathbb Q_{\ell}(1) ) \xrightarrow{\sim} T_{\ell}J_K\otimes \mathbb Q_{\ell}, \qquad H^1_{et}(C_K \times \bar{K}, \mathbb Q_{\ell})  \xrightarrow{\sim}  \textrm{Hom}(H^1_{et}(C_K \times \bar{K}, \mathbb Q_{\ell}(1)), \mathbb Q_{\ell}).\]
So we obtain an isomorphism
\begin{equation}\label{smooth} 
 H^1_{et}(C_K \times \bar{K}, \mathbb Q_{\ell}) \xrightarrow{\sim} \textrm{Hom}(T_{\ell}J_K\otimes \mathbb Q_{\ell}, \mathbb Q_{\ell})
\end{equation} 
of  $\textrm{Gal}(\bar{K}/K)
$-representations. Since $H^1_{et}(C_v \times \bar{\mathbb F}_q, \mathbb Q_{\ell})$ isomorphic to the subspace of $H^1_{et}(C_K \times \bar{K}, \mathbb Q_{\ell})$ of invariants under the inertia subgroup \cite[Lemma 1.2]{MR899399}, it follows from 
(\ref{smooth}) and (\ref{localfactor}) that $P_1(C_v, t)$ is $L_v(J, t)$.

Finally, one has the elementary identity  \footnote{This proposition, first stated on page 176 of \cite{MR528839}, has a typo which is corrected in its restatement on page 193.}  \cite[Proposition 3.3]{MR528839}; see also \cite[p.484]{MR2092767}: 
\[P_2(C_v,t) = \prod_{i\in I} (1-{(qt)}^{r_i}).\] 
This completes the proof.\end{proof} 

\subsection{Relating $\Sha(J/F)$ and $\textrm{Br}(X)$} For any arithmetic surface $X \to S$, let $\delta$ be the index of  $X_0$ over $F$ and $\alpha$ be the order of the (finite) cokernel of the natural map $\textrm{Pic}^0(X_0) \hookrightarrow J(F)$.  For any finite place $v$ of $S$, we put $\delta_v$ and $\delta'_v$ for the (local) index and period of $X\times F_v$ over the local field $F_v$. The following result is due to Geisser \cite[Theorem 1.1]{Geisser1}; there is also a recent proof by Flach-Siebel \cite{flachsiebel}.  
\begin{theorem}\label{geisser} Assume that ${\rm Br}(X)$ is finite. The following equality holds (up to powers of $2$):
\begin{equation}\label{konaaa} [{\rm Br}(X)] \alpha^2\delta^2 = [\Sha(J/F)] \prod_{v \in S} \delta'_v \delta_v.\end{equation}
\end{theorem}

\section{The proof of the main theorem in the general case}\label{Proof2} 

\subsection{Preliminary steps} We are in a position to prove Theorem \ref{Main theorem} in the general case. Let $\Sigma=\{v\in S~|~X_v~\textrm{is~not~smooth}\}$; let  $G_v$ denote the set of irreducible components of $X_v$.

\begin{lemma}\label{q12}  If
\[ Q_2(s) = \prod_{v \in \Sigma} \frac{(1-q_v^{1-s})}{\prod_{i\in G_v} (1-q_v^{{r_i}
(1-s)})},\] then 
 \begin{equation}\label{q12-ar}  Q_2^*(1) = \frac{1}{\prod_{v\in \Sigma} ( (\log q_v)^{\#G_v-1}~\cdot~\prod_{i\in G_v} r_i)} = \prod_{v\in \Sigma}\frac{c_v}{\Delta_{ar}(R_v)~.\delta_v~\cdot~\delta'_v} = P_{J,fin}~\cdot~\prod_{v\in \Sigma}\frac{1}{\Delta_{ar}(R_v)~.\delta_v~\cdot~\delta'_v}\end{equation}
\end{lemma} 
\begin{proof} The first equality is clear; the second uses Corollary \ref{fls}; the third follows from  (\ref{tamagawas}) using $c_v =1$ for $v \notin \Sigma$. \end{proof}
 
 \begin{proposition}\label{B-final} 
One has 
\[ \chi(B(1,1)) = \frac{P_{J, fin}}{P_{J}} = \frac{1}{P_{J,\infty}}=  \frac{|d_F|^{g/2}}{P_{J, \infty}(\eta)~\cdot~ \mathbb N_{F/{\mathbb Q}}(\mathfrak a_{\eta})}.\] 
\end{proposition} 
\begin{proof} One has the isomorphism of integral structures (its Euler characteristic is one)
\[H^1_B(J_{\mathbb C}, \mathbb Z(1))^{+}_{\mathbb C} \xrightarrow{\sim} H^1_B(X_{\mathbb C}, \mathbb Z(1))^{+}_{\mathbb C}.\]
 
 As the notion of N\'eron model is local on the base \cite[Proposition 4, page 13]{BLR}, Theorem \ref{raynaud1} shows that the kernel and cokernel of 
${\rm Lie}(\phi): H^1(X, \mathcal O) \to {\textrm{Lie}}(\mathcal J)$ are torsion $\mathcal O_F$-modules of the same length and  the kernel is exactly the torsion of $H^1(X,\mathcal O_X)$.   This implies that the Euler characteristic of 
\[H^1(X, \mathcal O)\otimes_{\mathbb Z} \mathbb C\xrightarrow{{\textrm{Lie}}{(\phi})} {\textrm{Lie}}(\mathcal J)\otimes_{\mathbb Z}\mathbb C\]
is one: $\chi({\rm Lie}(\phi)) =1$. Here the integral structures are $(H^1(X, \mathcal O), [H^1(X, \mathcal O)_{tor}])$ and $({\textrm{Lie}}(\mathcal J), 1)$.  
We obtain that the Euler characteristic of
\[ H^1_B(X_{\mathbb C}, \mathbb Z(1))^{+}_{\mathbb C} \xrightarrow{\gamma_X} H^1(X, \mathcal O)\otimes_{\mathbb Z} \mathbb C\]
is equal to that of 
\[\gamma_{\mathcal J}: H^1_B(J_{\mathbb C}, \mathbb Z(1))^{+}_{\mathbb C} \xrightarrow{\sim}H^1_B(X_{\mathbb C}, \mathbb Z(1))^{+}_{\mathbb C} \xrightarrow{\gamma_X} H^1(X, \mathcal O)\otimes_{\mathbb Z} \mathbb C\xrightarrow{{\textrm{Lie}}{(\phi})} {\textrm{Lie}}(\mathcal J)\otimes_{\mathbb Z}\mathbb C.\]

 Proposition \ref{yetanother} and Theorem \ref{allperiods} show
\begin{equation}\label{kone-illa2} \chi(\gamma_X) = \chi(\gamma_{\mathcal J}) = 
\frac{P_{J,\infty}(\eta)~\cdot~\mathbb N_{F/{\mathbb Q}}(\mathfrak a_{\eta})}{{\sqrt{|d_F|}}^g}.\end{equation}
Our convention for $B(1,1)$ is that $\textrm{Ker}(\gamma_M)$ is in even degree, say degree zero. Thus
\[\underset{\textrm{degree one}}{H^1_B(X_{\mathbb C}, \mathbb Z(1))^{+}_{\mathbb C}} \to \underset{\textrm{degree two}}{H^1(X, \mathcal O)\otimes \mathbb C}, \qquad det (B(1,1))  = \frac{{\sqrt{|d_F|}}^g}{P_{J,\infty}(\eta)~\cdot~\mathbb N_{F/{\mathbb Q}}(\mathfrak a_{\eta})} = \chi(B(1,1)). \]
\end{proof}

 \subsection{Global index}\label{indexooo}  Since $X$ is regular, the natural map $ \textrm{Pic}(X) \to \textrm{Pic}(X_0)$ is surjective: a natural section (as sets) is provided by sending a divisor on $X_0$ to its Zariski closure in $X$. So the index $\delta$ of $X_0$ over $F$ is the order of 
the cokernel of the composite map $\textrm{Pic}(X) \to \textrm{Pic}(X_0) \xrightarrow{d} \mathbb Z$.
\subsubsection{Calculation of $\chi(C(1))$} The torsion in $C(1)$ of (\ref{C1}) satisfies
 \begin{equation}\label{chiC1tors} \chi(C(1)_{tor}) =  \frac{w~\cdot~[Br(X)]} {[\textrm{Pic}(X)_{tor}]~\cdot~[\textrm{Pic}(X)_{tor}]}.\end{equation}

We can rewrite $C(1)_{\mathbb C}$ as
\begin{equation}\label{c1nue} 
0 \to  \underset{degree~zero}{\mathcal O_F^{\times}\otimes\mathbb C} \to \underset{degree~one}{\textrm{Coker}(\gamma_{M_0})} \to \underset{degree~two}{\textrm{Hom}(\textrm{Pic}(X), \mathbb C)} \to \underset{degree~three}{{\textrm{Pic}(X)}\otimes\mathbb C} \to \underset{degree~four}{\textrm{Ker}(\gamma_{M_2})} \to \underset{degree~five}{\textrm{Hom}(\mathcal O_F^{\times}, \mathbb C)} \to 0.
\end{equation} 
Using the exact sequence
\begin{equation}\label{globalindex} 0 \to \textrm{Pic}^0(X) \to \textrm{Pic}(X)  \xrightarrow{d} \mathbb Z \to \frac{\mathbb Z}{\delta\mathbb Z} \to 0,\end{equation}
the sequence $C(1)$ breaks up into 
\begin{itemize} 
\item the sequence  (\ref{conj-S-1}) with determinant $R$
\[0 \to  \underset{degree~zero}{\mathcal O_F^{\times}\otimes\mathbb C} \to \underset{degree~one}{\textrm{Coker}(\gamma_{M_0})} \to \underset{degree~two}{\textrm{Hom}(\mathbb Z, \mathbb C)} \to 0;\]
\item the sequence with determinant $\delta$
\[0 \to \textrm{Hom}(\mathbb Z, \mathbb C) \to \textrm{Hom}(\textrm{Pic}(X), \mathbb C)  \to \textrm{Hom}(\textrm{Pic}^0(X), \mathbb C) \to 0; \]
 \item and the sequence with determinant $det(\psi)$
 \[h: \underset{degree~two}{\textrm{Hom}(\textrm{Pic}^0(X), \mathbb C)} \to \underset{degree~three}{{\textrm{Pic}^0(X)}\otimes\mathbb C} \]
 where $h$ is the Arakelov intersection pairing \cite[Conjecture 2.3.7]{SL2020} as in (\ref{Cr});
 \item the sequence with determinant $\delta$
 \[0 \to \textrm{Pic}^0(X)\otimes \mathbb C \to \textrm{Pic}(X)\otimes \mathbb C  \xrightarrow{d} \mathbb C \to 0;\]
 \item the sequence (\ref{conj-S-zero}) with determinant $R$
\[ 0 \to \mathbb C \to \underset{degree~four}{\textrm{Ker}(\gamma_{M_2})} \to \underset{degree~five}{\textrm{Hom}(\mathcal O_F^{\times}, \mathbb Z)\otimes \mathbb C} \to 0.\]
 \end{itemize} 
 It follows from (\ref{globalindex}) that, in (\ref{c1nue}), the image of the lattice in degree three has index $\delta$ in the lattice in degree four; dually, the image of the lattice in degree one has index $\delta$ in the lattice in degree three. Thus, 
\[det(C(1)) =  \frac{R~\cdot~R}{\delta~\cdot~\delta~\cdot~det(h)}.\] 
So  $ \chi(C(1)) $ is given as  
\begin{equation}\label{detc1}= \frac{det(C(1))}{\chi(C(1)_{tor})} =  \frac{R^2}{\delta^2~\cdot~\Delta_{ar}(\textrm{Pic}^0(X))~\cdot~[\textrm{Pic}^0(X)_{tor}]^2}~\cdot~\frac{[\textrm{Pic}^0(X)_{tor}]^2}{w~\cdot~[\textrm{Br}(X)]} = \frac{R^2}{\delta^2~\cdot~\Delta_{ar}(\textrm{Pic}^0(X))~\cdot~w~\cdot~[\textrm{Br}(X)]},\end{equation} 
 where, as in \S \ref{pairings}, 
\begin{equation}\label{pic-ar} \Delta_{ar}(\textrm{Pic}^0(X)) = \frac{det(h)}{[\textrm{Pic}^0(X)_{tor}]^2}.\end{equation}
Our next task is to calculate $\Delta_{ar}(\textrm{Pic}^0(X))$ and relate it to the N\'eron-Tate pairing (\ref{jf-ar}) on $J(F)$. 
\subsubsection{Calculation of $\Delta_{ar}(\textrm{Pic}^0(X))$.} This is based on localization sequences on $X$ and $S$. 

Let $U = S - \Sigma$. So the map $X_U = \pi^{-1}(U) \to U$ is smooth. For any finite $\Sigma' \subset S$ containing $\Sigma$, we put $U' = S - \Sigma'$ and $X_{U'} = X - \pi^{-1}(U')$. 
\begin{lemma}\label{lastlemma}  (i) The maps 
\[\textrm{Pic}(S) \to  \textrm{Pic}(X), \qquad \textrm{Pic}(U') \to \textrm{Pic}(X_{U'})\]
are injective. 

(ii)   There is an exact sequence
\begin{equation}\label{anotherses} 0 \to  \oplus_{v\in \Sigma} \frac{ \mathbb Z^{G_v}}{\mathbb Z} \to \frac{\textrm{Pic}^0 (X)}{\textrm{Pic}(S)} \to \textrm{Pic}^0 ( X_0) \to 0.\end{equation}
\end{lemma} 
\begin{proof}  (i) 
 From the Leray spectral sequence for $\pi: X\to S$ and the \'etale sheaf $\mathbb G_m$ on $X$, we get the exact sequence
 \[ 0 \to H^1(S, \pi_* \mathbb G_m) \to H^1(X, \mathbb G_m) \to H^0(S, R^1\pi_*\mathbb G
 _m) \to \textrm{Br}(S).\] 
Now use that $\pi_*\mathbb G_m$ is the sheaf $\mathbb G_m$ on $S$. 
This provides the injectivity of the first map. A similar argument provides the injectivity of the second.

(ii) We can compare the localization sequences for $X_{U'} \subset X$ and $U' \subset S$

\begin{tikzcd}
0 \arrow[r,] &\Gamma({X}, \mathbb G_m)   \arrow[r, ] \isoarrow{d} &\Gamma({U'}, \mathbb G_m) \arrow[r,] \isoarrow{d} &\underset{v\in Z'}{\oplus} \mathbb Z \arrow[r, ] \arrow[d, hook, red] & \textrm{Pic}(S)   \arrow[r, ] \arrow[d, hook, red] &  \textrm{Pic}(U') \arrow[r,] \arrow[d, hook, red]& 0\\
0 \arrow[r,] & \Gamma({X}, \mathbb G_m)  \arrow[r,] &\Gamma(X_{U'}, \mathbb G_m) \arrow[r,]  &\underset{v\in Z'}{\oplus} \mathbb Z^{G_v} \arrow[r,]  & \textrm{Pic}(X)   
\arrow[r,]   &  \textrm{Pic}(X_{U'}) \arrow[r,]  & 0
\end{tikzcd}

Note that $S$ and $X$ are regular and for any regular scheme $Y$, one has an isomorphism $\textrm{Cl}(Y) = \textrm{Pic}(Y)$ between the class group and the Picard group. 
The natural map between the  localization sequences is injective on all terms and, by assumption, is an  isomorphism on the first and second terms. This provides the exact sequence
\[0 \to  \oplus_{v\in \Sigma'} \frac{ \mathbb Z^{G_v}}{\mathbb Z} \to \frac{\textrm{Pic}(X)}{\textrm{Pic}(S)} \to \frac{\textrm{Pic}(X_{U'})}{\textrm{Pic}(U')} \to 0.\]
In particular, we get this sequence for $\Sigma$ and $U$. Using (\ref{globalindex}), we obtain the exact sequence
\[0 \to  \oplus_{v\in \Sigma} \frac{ \mathbb Z^{G_v}}{\mathbb Z} \to \frac{\textrm{Pic}^0(X)}{\textrm{Pic}(S)} \to \frac{\textrm{Pic}^0(X_{U})}{\textrm{Pic}(U)} \to 0.\]
By assumption, $X_v$ is geometrically irreducible for any $v \notin \Sigma$. 
So for any $U' = S - \Sigma'$ with  $U' \subset U$, the induced maps 
\[\frac{\textrm{Pic}^0(X_U)}{\textrm{Pic}(U)}\to \frac{\textrm{Pic}^0(X_{U'})}{\textrm{Pic}(U')}\]
are isomorphisms. Taking the limit over $\Sigma'$ gives us an exact sequence
\[0 \to  \oplus_{v\in \Sigma} \frac{ \mathbb Z^{G_v}}{\mathbb Z} \to \frac{\textrm{Pic}^0 (X)}{\textrm{Pic}(S)} \to \textrm{Pic}^0(X_0) \to 0.\] This proves the lemma.\end{proof} 

For any $v \in S$, recall  $\Delta_{ar}(R_v)$ from Corollary \ref{fls} where 
 \[R_v =\frac{\mathbb Z^{G_v}}{\mathbb Z},\]
and, as before, $G_v$ is the set of irreducible components of $X_v$.  Let us define (see \S \ref{pairings}) 
 \begin{equation}\label{jf-ar} \Delta_{NT}(J(F)) = \frac{\Theta_{NT}(J)}{ [ J(F)_{tor}]~\cdot~[J(F)_{tor}]} \end{equation} using the N\'eron-Tate pairing (\ref{Nerontate}) on $J(F)$, analogous to $\Delta_{ar}(\textrm{Pic}^0(X))$ from (\ref{pic-ar}).  

 \begin{proposition}\label{arakelovpairing}  If $\alpha$ is the order of the cokernel of the natural map $\textrm{Pic}^0(X_0) \hookrightarrow J(F)$, 
then one has 
 \begin{equation}\label{kone-illa3} \Delta_{ar}(\textrm{Pic}^0(X))=  \pm \frac{\alpha^2}{h^2}~\cdot~\Delta_{NT}(J(F))~\cdot~\prod_{v\in\Sigma} \Delta_{ar}(R_v)\end{equation}
 \end{proposition} 
\begin{proof} By \cite[Proposition 3.3]{MR778087}, the Arakelov intersection pairing $h$ is negative-definite on $\textrm{Pic}^0(X)\otimes \mathbb Q$, there exists a map $\kappa: \textrm{Pic}^0(X_0) \to \textrm{Pic}^0(X)$ such that 
\[ (y,y') \mapsto h(\kappa(y), \kappa(y') )\]
gives the intersection pairing on $\textrm{Pic}^0(X_0)$ which, by Faltings-Hriljac \cite{MR740897}, \cite[Theorem 3.1]{MR778087}, is the negative of the N\'eron-Tate pairing (\ref{Nerontate}) on $J(F)$. Note that this means 
\begin{equation}\label{twodeltas} 
\Delta_{NT}(J(F)) = \pm \Delta_{ar}(J(F)).
\end{equation}
So the sequence  (\ref{anotherses}) splits over $\mathbb Q$ as an orthogonal direct sum with respect to the Arakelov intersection pairing:  
\begin{equation}\label{orthop} (\textrm{Pic}^0 ( X_0)\otimes\mathbb Q)~\oplus~(\underset{v\in \Sigma}{\oplus} \frac{ \mathbb Q^{G_v}}{\mathbb Q}) \cong  \textrm{Pic}^0 (X)\otimes\mathbb Q.\end{equation}
The map $\kappa$ is defined as follows: given any element $y$ of $\textrm{Pic}^0 ( X_0)$, consider  its Zariski closure $\bar{y}$ in $X$. As the intersection pairing is negative-definite  \cite[\S 9.1, Theorem 1.23]{MR1917232} on $R_v$, the linear mapping $R_v \to \mathbb Z$ defined by $z \mapsto z.\bar{y}$ is represented by a unique element $\kappa_v(y) \in R_v\otimes\mathbb Q$. Clearly, the element
\[\kappa(y)  = \bar{y} - \sum_{v\in \Sigma} \kappa_v(y) \]
is orthogonal to $\oplus_{v\in \Sigma}R_v \subset \textrm{Pic}^0(X)$; so the assignment $y \mapsto \kappa(y)$ provides (\ref{orthop}).

Now (\ref{deltann}) shows 
\[ \Delta_{ar}(\frac{\textrm{Pic}^0(X)}{\textrm{Pic}(S)})=  \Delta_{ar}(\textrm{Pic}^0(X_0 ))~\cdot~\prod_{v\in\Sigma} \Delta_{ar}(R_v).\]
Using $h = [\textrm{Pic}(S)]$, this becomes  
\begin{equation}\label{hpic}\Delta_{ar}(\textrm{Pic}^0(X))~\cdot~ {h^2} = \Delta_{ar}(\frac{\textrm{Pic}^0(X)}{\textrm{Pic}(S)})=  \Delta_{ar}(\textrm{Pic}^0(X_0 ))~\cdot~\prod_{v\in\Sigma} \Delta_{ar}(R_v).\end{equation}
As $\textrm{Pic}^0(X_0) \hookrightarrow J(F)$ is a subgroup of index $\alpha$, we see $ \Delta_{ar}(\textrm{Pic}^0(X_0 )) = \alpha^2 \Delta_{ar}(J(F)) = \pm \alpha^2\cdot \Delta_{NT}(J(F))$.\end{proof}

\subsection{The zeta function of $X$} Our first step is to rewrite the zeta function $\zeta(X,s)$
\[ \zeta(X,s) = \prod_{v \in S}\zeta(X_v, q_v^{-s}), \qquad \zeta(X_v, t) = \frac{P_1(X_v,t)}{P_0(X_v, t) P_2(X_v, t)}.\] 

\begin{proposition}\label{newprop} One has 
\begin{equation}\label{kone-illa} \zeta(X,s)  = \frac{\zeta(S,s)~\cdot~\zeta(S, s-1)}{L(J,s)}~\cdot~Q_2(s), \qquad \zeta^*(X,1) = \frac{\zeta^*(S,1)~\cdot~\zeta^*(S, 0)}{L^*(J,1)}~\cdot~Q_2^*(1).\end{equation}
\end{proposition}  
\begin{proof} 
By Lemma \ref{cohflat}, we see that 

(i) the factors $P_0$ combine to give $\zeta(S,s)$. 

(ii) the factors $P_1$ combine to give $L(J,s)$ .

(iii)  $P_2(X_v,t)$ is the expected factor $(1-q_vt)$ for $v\notin \Sigma$; for $v \in \Sigma$, it is given in Lemma \ref{cohflat}.  So
\[\prod_{v\in S}\frac{1}{P_2(X_v,q_v^{-s}) } = Q_2(s)~\cdot~\prod_{v\in S} \frac{1}{1-q_v^{(1-s)}} = Q_2(s)~\cdot~\zeta(S, s-1).\]
\end{proof}

\subsection{Completion of the proof of Theorem \ref{Main theorem}} We prove the main result of this paper. 
\begin{proof} (of Theorem \ref{Main theorem}) We can restate Conjecture \ref{bsd} for $J/F$ using (\ref{jf-ar}) as 
\begin{equation}\label{bsd-jf-ar} \frac{1}{L^*(J,1)} = \frac{|d_F|^{g/2}}{P_{J, fin}~\cdot~P_{J, \infty}(\eta)~\cdot~ \mathbb N_{F/{\mathbb Q}}(\mathfrak a_{\eta})}~\cdot~  \frac{1}{\Delta_{NT}(J(F))~\cdot~[\Sha(J/F)]}.\end{equation} 

From Proposition \ref{newprop}, we see that Conjecture \ref{bsd} for $J/F$
 is equivalent to the equality
\begin{equation}\label{g=0} 
\begin{split} 
\zeta^*(X,1) & = \zeta^*(S,1)~\cdot~\zeta^*(S, 0)~\cdot~\frac{1}{L^*(J,1)}~\cdot~Q_2^*(1) \\
& \overset{(\ref{cnf1})}{=} \frac{2^{r_1} (2\pi)^{r_2}hR}{w \sqrt{|d_F|}} ~\cdot~{\zeta^*(S, 0)}~\cdot~\frac{1}{L^*(J,1)}~\cdot~Q_2^*(1) \\ 
& \overset{(\ref{cnf0})}{=} \frac{2^{r_1} (2\pi)^{r_2}hR}{w \sqrt{|d_F|}}~\cdot~\frac{hR}{w}~\cdot~\frac{1}{L^*(J,1)}~\cdot~Q_2^*(1) \\
&\overset{(\ref{bsd-jf-ar})}{=} \frac{2^{r_1} (2\pi)^{r_2}hR}{w \sqrt{|d_F|}}~\cdot~\frac{hR}{w}~\cdot~ \frac{|d_F|^{g/2}}{P_{J, fin}~\cdot~P_{J, \infty}(\eta)~\cdot~ \mathbb N_{F/{\mathbb Q}}(\mathfrak a_{\eta})}~\cdot~  \frac{1}{\Delta_{ar}(J(F))~\cdot~[\Sha(J/F)]}~\cdot~Q_2^*(1) \\
& \overset{(\ref{q12-ar})}{=} \frac{2^{r_1} (2\pi)^{r_2}h^2R^2\cdot |d_F|^{g/2} \cdot P_{J,fin} }{w^2 \sqrt{|d_F|} \cdot P_{J, fin} \cdot P_{J, \infty}(\eta) \cdot \mathbb N_{F/{\mathbb Q}}(\mathfrak a_{\eta})  \cdot \Delta_{ar}(J(F)) \cdot [\Sha(J/F)]} \cdot \prod_{v\in \Sigma} \left( \frac{1}{\Delta_{ar}(R_v) \cdot \delta_v \cdot\delta'_v}\right)\\
&= \frac{2^{r_1} (2\pi)^{r_2}}{w \sqrt{|d_F|}} \cdot \frac{h^2 \cdot R^2}{w} \cdot \frac{|d_F|^{g/2}}{P_{J, \infty}(\eta)~\cdot~ \mathbb N_{F/{\mathbb Q}}(\mathfrak a_{\eta})} \cdot   
\frac{1}{\Delta_{ar}(J(F)) \cdot [\Sha(J/F)]} \cdot \prod_{v\in \Sigma} \left(\frac{1}{\Delta_{ar}(R_v) \cdot \delta_v \cdot~ \delta'_v}\right)
\end{split}
\end{equation}

 Conjecture \ref{main conjecture} states that
 \[ \zeta^*(X,1) = \chi(X,1) = \frac{\chi_{A,C}(X,1)}{\chi_B(X,1)} = \chi_{A}(X,1)~\cdot~\chi_{A'}(X,1)~\cdot~\chi(C(1))~\cdot~\frac{\chi(B(1,1))}{\chi(B(0,1)) \chi(B(2,1))}.\] 
 Using (\ref{bzero1}), (\ref{btwo1}),  (\ref{chiaa'}) and Proposition \ref{B-final} 
\[ \chi(B(0,1)) =  \frac{\sqrt{|d_F|}}{(2\pi )^{r_2}},\quad  \chi(B(1,1)) = \frac{|d_F|^{g/2}}{P_{J, \infty}(\eta)~\cdot~ \mathbb N_{F/{\mathbb Q}}(\mathfrak a_{\eta})},\quad  \chi(B(2,1)) =1= \chi_A(X,1), \quad \chi_{A'}(X,1) = \frac{1}{w},\]
we have 
\begin{equation}\label{zeta1fullformula} 
\begin{split}
\chi(X,1) &=  \frac{\chi_{A,C}(X,1)}{ \chi_B(X,1)} \\
&=  \chi_{A}(X,1)~\cdot~\chi_{A'}(X,1)~\cdot~\chi(C(1))~\cdot~\frac{\chi(B(1,1))}{\chi(B(0,1)) \chi(B(2,1))} \\
 &= 1.~\frac{1}{w}~\cdot~\chi(C(1))~\cdot~\frac{(2\pi )^{r_2}}{\sqrt{|d_F|}} ~\cdot~\frac{|d_F|^{g/2}}{P_{J, \infty}(\eta)~\cdot~ \mathbb N_{F/{\mathbb Q}}(\mathfrak a_{\eta})}.\\
 &= \chi(C(1))~\cdot~\frac{(2\pi )^{r_2}}{w\sqrt{|d_F|}} ~\cdot~\frac{|d_F|^{g/2}}{P_{J, \infty}(\eta)~\cdot~ \mathbb N_{F/{\mathbb Q}}(\mathfrak a_{\eta})}\\
 & \overset{(\ref{detc1})}{=} \frac{R^2}{\delta^2~\cdot~\Delta_{ar}(\textrm{Pic}^0(X))~\cdot~w~\cdot~[\textrm{Br}(X)]}~\cdot~\frac{(2\pi )^{r_2}}{w\sqrt{|d_F|}} ~\cdot~\frac{|d_F|^{g/2}}{P_{J, \infty}(\eta)~\cdot~ \mathbb N_{F/{\mathbb Q}}(\mathfrak a_{\eta})}
 \end{split}
 \end{equation}
 To prove Theorem \ref{Main theorem}, we need to show that $\chi(X,1)$
is equal to 
\[ \frac{2^{r_1} (2\pi)^{r_2}}{w \sqrt{|d_F|}} \frac{h^2~\cdot~R^2}{w}~\cdot~ \frac{|d_F|^{g/2}}{P_{J, \infty}(\eta)~\cdot~ \mathbb N_{F/{\mathbb Q}}(\mathfrak a_{\eta})}~\cdot~  
\frac{1}{\Delta_{ar}(J(F))~\cdot~[\Sha(J/F)]}~\cdot~\prod_{v\in \Sigma}~\frac{1}{\Delta_{ar}(R_v)~\cdot~\delta_v~\cdot~\delta'_v}.\]

Thus, we see that Theorem \ref{Main theorem} follows if we prove that (neglecting powers of two)
\begin{equation}
\begin{split}
\chi(C(1)) &\overset{?}{=} \frac{h^2~\cdot~R^2}{w}~\cdot~\frac{1}{\Delta_{ar}(J(F))~\cdot~[\Sha(J/F)]}~\cdot~\prod_{v\in \Sigma}~\frac{1}{\Delta_{ar}(R_v)~\cdot~\delta_v~\cdot~\delta'_v}\\
& \overset{?}{=} \frac{R^2}{w}~\cdot~\frac{h^2}{\Delta_{ar}(J(F))}~\cdot~\prod_{v\in \Sigma}~\frac{1}{\Delta_{ar}(R_v)} ~\cdot~\frac{1}{[\Sha(J/F)]}~\cdot~\prod_{v\in \Sigma}~\frac{1}{\delta_v~\cdot~\delta'_v}
\end{split}
\end{equation}
By (\ref{detc1}), we have
\[\chi(C(1)) = \frac{R^2}{\delta^2~\cdot~\Delta_{ar}(\textrm{Pic}^0(X))~\cdot~w~\cdot~[\textrm{Br}(X)]}.\] 
Thus, to prove Theorem \ref{Main theorem}, we only need to check if
\[  \frac{1}{\delta^2~\cdot~\Delta_{ar}(\textrm{Pic}^0(X))~\cdot~[\textrm{Br}(X)]} \overset{?}{=} \frac{h^2}{\Delta_{ar}(J(F))}~\cdot~\prod_{v\in \Sigma}~\frac{1}{\Delta_{ar}(R_v)} ~\cdot~\frac{1}{[\Sha(J/F)]}~\cdot~\prod_{v\in \Sigma}~\frac{1}{\delta_v~\cdot~\delta'_v}.\]
We can verify this readily using Proposition \ref{arakelovpairing} which says
\begin{equation}\label{g=zero} \Delta_{ar}(\textrm{Pic}^0(X))=  \frac{\alpha^2~\cdot~\Delta_{ar}(J(F))~\cdot~\prod_{v\in\Sigma} \Delta_{ar}(R_v)}{h^2}\end{equation}
and Theorem \ref{geisser} which says (note $\delta_v =1 = \delta'_v$ if $v \notin \Sigma$)
\begin{equation}\label{geisser2} [{\rm Br}(X)] \alpha^2\delta^2 = [\Sha(J/F)] \prod_{v \in S} \delta'_v \delta_v.\end{equation}
This completes the proof of Theorem \ref{Main theorem} in the general case.\end{proof} 

\begin{remark}\label{giszero} Consider the case $g=0$. As $J$ is trivial, Theorem \ref{Main theorem} says Conjecture \ref{main conjecture} is true. Let us see directly that this is true (up to powers of two): (\ref{g=0}) shows that
\[\zeta^*(X,1) = \frac{2^{r_1} (2\pi)^{r_2}}{w \sqrt{|d_F|}}~\cdot~\frac{h^2~\cdot~R^2}{w} ~\cdot~\prod_{v\in \Sigma}~\frac{1}{\Delta_{ar}(R_v)~\cdot~\delta_v~\cdot~\delta'_v},\]
but  (\ref{zeta1fullformula}) and (\ref{g=zero}) show
\[\chi(X,1) =   \frac{R^2}{\delta^2~\cdot~\Delta_{ar}(\textrm{Pic}^0(X))~\cdot~w~\cdot~[\textrm{Br}(X)]}~\cdot~\frac{(2\pi )^{r_2}}{w\sqrt{|d_F|}} = \frac{h^2}{\alpha^2~\cdot~\prod_{v\in \Sigma}~{\Delta_{ar}(R_v)}}~\cdot~\frac{R^2}{\delta^2~\cdot~w~\cdot~[\textrm{Br}(X)]}~\cdot~\frac{(2\pi )^{r_2}}{w\sqrt{|d_F|}}.\]
As $\alpha$, $\delta_v$,  $\delta$ and $\delta'_v$ divide two and (\ref{geisser2}) shows that $[\textrm{Br}(X)]$ is a power of two, it follows that the equality $\zeta^*(X,1) = \chi(X,1)$ is valid up to a finite power of two.\qed
\end{remark} 

\subsection{Comparison with the Artin-Tate conjecture}\label{ATooo} Theorem \ref{Main theorem} is the analogue for arithmetic surfaces of Conjecture (d) \cite[p.~427]{Tate1966} in the spirit of \cite{SL1983}.  More precisely, one has the
\begin{theorem} 
 Conjecture \ref{main conjecture} for $X$ is equivalent to the following identity (up to a finite power of two)
 \begin{equation}\label{atate} \frac{\zeta^*(S,0)\cdot \zeta^*(S,1)}{\zeta^*(X,1)} = \chi(\gamma_X)\cdot [\textrm{Br}(X)]\cdot \Delta_{ar}(\frac{\textrm{Pic}^0(X)}{\textrm{Pic}(S)}) \cdot \delta^2.
 \end{equation} 
\end{theorem} 
Here $\delta$ is the global index of $X$ as in \S \ref{indexooo}. The term $\chi(\gamma_X)$ is an {\it archimedean period}: it is the Euler characteristic (as in \S \ref{integral1}) of  
\begin{equation*} \gamma_X: H^1_B(X_{\mathbb C}, \mathbb Z(1))^{+}_{\mathbb C} \xrightarrow{\sim} H^1(X, \mathcal O)\otimes_{\mathbb Z} \mathbb C.
\end{equation*}  with respect to the integral structures provided by $H^1_B(X_{\mathbb C}, \mathbb Z(1))^{+}$ and the abelian group underlying $H^1(X, \mathcal O)$.
\begin{proof} We can start with the left hand side and rewrite it as follows:
\begin{equation*} 
\begin{split}
 \frac{\zeta^*(S,0)\cdot \zeta^*(S,1)}{\zeta^*(X,1)} &\overset{}{=} \frac{\zeta^*(S,0)\cdot \zeta^*(S,1)}{\chi(X,1)} = \frac{\zeta^*(S,0)\cdot \zeta^*(S,1) \cdot\chi_B(X,1)}{\chi_A(X,1) \cdot \chi_{A'}(X,1)\cdot \chi(C(1))}\\
 &= \frac{\zeta^*(S,0)\cdot \zeta^*(S,1) \cdot\chi(\gamma_X)}{1\cdot \frac{1}{w}\cdot \chi(C(1)) \cdot \frac{{(2\pi)}^{r_2}}{\sqrt{|d_F|}}} \quad \textrm{by}~(\ref{bzero1}), (\ref{btwo1}),  (\ref{chiaa'}) \\
 & {=}  \frac{\chi(\gamma_X) \cdot \zeta^*(S,0)\cdot \zeta^*(S,1) \cdot w^2 \cdot \delta^2\cdot \Delta_{ar}(\textrm{Pic}^0(X))\cdot [\textrm{Br}(X)] \cdot \sqrt{|d_F|}}{R^2\cdot (2\pi)^{r_2} } 
 \quad \textrm{by}~(\ref{detc1})\\
 &= \chi(\gamma_X) \cdot \delta^2\cdot h^2\cdot \Delta_{ar}(\textrm{Pic}^0(X))\cdot [\textrm{Br}(X)] \cdot 2^{r_1} \quad \textrm{by}~(\ref{cnf1}), (\ref{cnf0}) \\
 &{=}  \chi(\gamma_X) \cdot \delta^2\cdot \Delta_{ar}\left(\frac{\textrm{Pic}^0(X)}{\textrm{Pic}(S)} \right)\cdot [\textrm{Br}(X)] \quad \textrm{by}~(\ref{hpic}).
 \end{split}
 \end{equation*} 
 This proves the stated equivalence. \end{proof} 
 \begin{remark}\label{suzuki} (T. Suzuki) The identity (\ref{atate}) is an analogue for arithmetic surfaces of the Artin-Tate conjecture \cite[Conjecture C]{Tate1966}. Let us show this by rewriting the Artin-Tate conjecture for function fields \cite[Conjecture 2]{LRS}; we shall use the notation of \cite{LRS} from now on.  
 
 So now let $S$ be a smooth proper curve over a finite field $F_q$.  And let $X$ be a smooth proper surface with a flat proper morphism $\pi:X \to S$ whose generic fiber  is a smooth geometrically connected curve $X_0$ over Spec~$\mathbb F_q(S)$. By combining (2), (12), and (13) of \cite{LRS}, we obtain
 \[\frac{\zeta(S,1-s) \cdot \zeta(S, s)}{\zeta(X,{s})} = \frac{P_2(X, q^{-s})}{P_1(B, q^{-s}) P_1(B, q^{1-s}) (1 - q^{1-s})^2}.\]
 Hence the Artin-Tate conjecture  \cite[Conjecture 2]{LRS} is equivalent to the identity
 \[\frac{\zeta^*(S,0)\cdot \zeta^*(S, 1)}{\zeta^*(X,1} = \frac{P_2^*(X, \frac{1}{q}) \cdot q^{\mathrm{dim}~B}}{ [B(\mathbb F_q)]^2\cdot (\log q)^2}  
 = \frac{[\textrm{Br}(X)]\cdot \Delta_{ar}(NS(X))\cdot q^{\chi(S, R^1\pi_*\mathcal O_X)}}{ [B(\mathbb F_q)]^2 \cdot (\log q)^2} \]
 by (14), (3) and (4) of \cite{LRS}. On the other hand, by \cite[Proposition 13 and Corollary 5 (ii)]{LRS}, we have
 \[\frac{\Delta_{ar}(NS(X))}{(\log q)^2} = \Delta_{ar}\left(\frac{NS(X)_0}{\pi^*NS(S)}\right) \cdot \delta^2 =  [B(\mathbb F_q)]^2\cdot \Delta_{ar} \left(\frac{\textrm{Pic}(X)_0}{\pi^*\textrm{Pic}(S)}\right) \cdot \delta^2.\]
Thus the Artin-Tate conjecture for $X$ is equivalent to the identity
 \[\frac{\zeta^*(S,0) \cdot \zeta^*(S,1)}{\zeta^*(X,1)} = q^{\chi(S, R^1\pi_*\mathcal O_X)}\cdot [\textrm{Br}(X)]\cdot \delta^2\cdot \Delta_{ar}\left(\frac{\textrm{Pic}(X)_0}{\pi^*\textrm{Pic}(S)}\right). \] 
 The factor $q^{\chi(S, R^1\pi_*\mathcal O_X)}$  plays the role of $\chi(\gamma_X)$ in (\ref{atate}).\qed
 \end{remark} 
Competing interests: The author(s) declare none


\end{document}